\documentclass[journal,final]{IEEEtran}

\usepackage{bm,amsmath,amssymb,graphicx,amsfonts,xspace}
\usepackage{psfrag,enumerate,pdfsync}
\usepackage[all,2cell,dvips]{xy} 
\newcommand{\comments}[1]{}

\usepackage{tikz}
\usetikzlibrary{arrows}

\usepackage{xcolor,framed}
\colorlet{shadecolor}{yellow}

\DeclareMathOperator{\diag}{diag}

\newcommand{\SO}{\mathsf{SO}}

\newcommand{\OO}{{\pmb{\boldsymbol \omega}}}
\newcommand{\bb}{{\boldsymbol \beta}}
\newcommand{\one}{{\boldsymbol 1}}

\def\qed{\hfill {$\square$}}

\newcommand{\cB}{{\cal B}}
\newcommand{\cI}{{\cal I}}

\renewcommand{\Re}{\mathbb{R}}
\newcommand{\cX}{{\cal X}}

\newcommand{\sR}{{\sf R}}
\newcommand{\sX}{{\sf X}}
\newcommand{\sO}{{\sf \Omega}}
\newcommand{\cN}{{\cal N}}

\newcommand{\trans}{{}^\top}

\newcommand{\wdes}{\pmb{\mathbf{f}}}
\newcommand{\gdes}{\pmb{\mathbf{g}}}
\newcommand{\hdes}{\pmb{\mathbf{h}}}
\newcommand{\kdes}{\pmb{\mathbf{k}}}
\newcommand{\lam}{\rho}
\newcommand{\kq}{k_1}
\newcommand{\ktilde}{k_2}
\newcommand{\w}{f}
\newcommand{\prj}{\theta}
\newcommand{\Chi}{X}
\newcommand{\Je}{\mathbb{J}}
\newcommand{\rot}{{\sf rot}}
\newcommand{\tran}{{\sf tran}}

\newtheorem{theorem}{Theorem}
\newtheorem{lemma}{Lemma}

\newtheorem{definition}{Definition}
\newtheorem{remark}{Remark}

\long\def\symbolfootnote[#1]#2{\begingroup%
\def\thefootnote{\fnsymbol{footnote}}\footnote[#1]{#2}\endgroup}

\title{Local and Distributed Rendezvous of Underactuated Rigid Bodies} %in
\author{Ashton Roza, Manfredi Maggiore, Luca Scardovi%%
\thanks{This research was
  supported by the National Sciences and Engineering Research Council
  of Canada.}%%
\thanks{The authors are with the Department of Electrical and
  Computer Engineering, University of Toronto, 10 King's College Road,
  Toronto, ON, M5S 3G4, Canada.  {\tt\footnotesize
    ashton.roza@mail.utoronto.ca, maggiore@ece.utoronto.ca, scardovi@scg.utoronto.ca}}}
\begin{document}
\maketitle

\begin{abstract}
This paper solves the rendezvous problem for a network of underactuated rigid bodies such as quadrotor helicopters. A control
strategy is presented that makes the centres of mass of the vehicles
converge to an arbitrarily small neighborhood of one another. The
convergence is global, and each vehicle can compute its own control
input using only an on-board camera and a three-axis rate
gyroscope. No global positioning system is required, nor any
information about the vehicles' attitudes.
\end{abstract}

\section{Introduction}\label{sec:intro}

\begin{figure}[t]
\centerline{\includegraphics[width=0.4\textwidth]{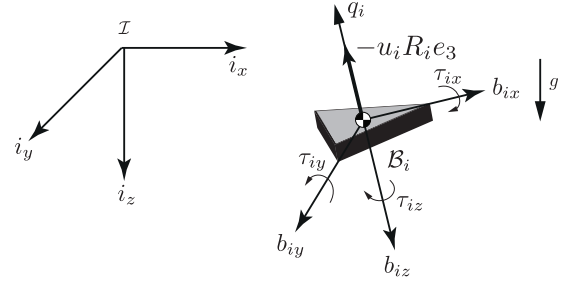}}
\caption{Vehicle class under consideration.} \label{fig:vehicle}
\end{figure}

Consider a network of flying robots, each propelled by a
thrust vector and endowed with an actuation mechanism producing
torques about three orthogonal body axes ---see
Figure~\ref{fig:vehicle}. With six degrees-of-freedom and four actuators, each robot is underactuated with degree of underactuation two. A quadrotor helicopter is an example of
such a robot. Suppose each robot mounts a camera and an inertial
measurement unit (IMU) that includes a three-axis rate-gyroscope, so
that the robot is able to measure, in the coordinates of its own
frame, the relative displacements and velocities of nearby vehicles,
and its own angular velocity.  The rendezvous control problem is to
get the robots to move to a common location using only the above on-board
sensors. To this day, this problem is open. This paper presents the
first solution.

Consider now $n \geq 2$ robots. The {\em
  rendezvous control problem} investigated in this paper is to find
feedback laws making the relative
distances and velocities become arbitrarily
small for all $i,j \in \{1,\dots,n\}$, and for arbitrary initial
conditions of all robots. Crucial in the problem statement is the
requirement on sensing.  If robot $i$ can sense robot $j$, then robot
$i$ can sense the relative position and velocity of robot $j$ in its
own local frame. Robot~$i$ can also measure its own angular
velocity in the coordinates of its body frame. Robot $i$ can neither access its own inertial position
and velocity, nor its own attitude.  A
feedback law satisfying the above sensing requirements is referred to
as being {\em local and distributed}.  In this paper, the set of
vehicles that robot $i$ can sense is assumed to
be constant. This assumption is questionable in practice, but is made
to render the problem mathematically treatable. The rendezvous problem
with distance-dependent neighbors remains a challenging open problem
for much simpler classes of robot models, such as double-integrators.
\begin{figure*}[t]
\centerline{\includegraphics[width=0.6\textwidth]{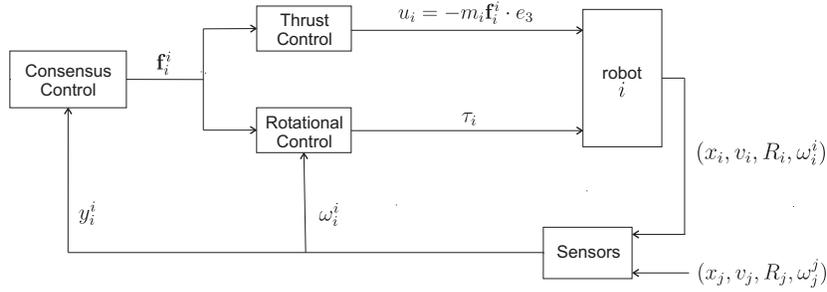}}
\caption{Block diagram of the rendezvous control system for robot
  $i$. The outer loop assigns a desired thrust vector
  $\wdes_i(y_i^i)$. The inner loop thrust control uses
  $\wdes_i(y_i^i)$ to assign the vehicle input $u_i$ while the
  rotational control uses $\wdes_i(y_i^i)$ to assign the torque input
  $\tau_i$. The vector $y_i^i$ contains the relative displacements and
  velocities of vehicles that robot $i$ can sense, measured in the body
  frame of robot $i$.}
\label{fig:block_diagram}
\end{figure*}
The block diagram of the proposed controller is depicted in
Figure~\ref{fig:block_diagram}. There are two nested loops. The outer
loop treats each robot as a point-mass driven by a force input, and
produces a double-integrator consensus controller which
becomes a reference input for the inner loop. The inner loop assigns
local and distributed feedbacks for the robots. More intuition
is provided in Section~\ref{sec:soln_pcp_2stage}.

Besides having a simple expression making its real-time implementation
feasible, the proposed controller meets the sensing requirements of
the rendezvous control problem. In particular, it does not require any
knowledge of the robots' absolute positions and velocities, or of
their attitudes. It does not even require sensing of the relative
attitudes. Finally, the controller does not require any communication
among robots.

Our main result, Theorem~\ref{thm:2step}, states that the proposed controller does indeed solve the rendezvous control problem, and in so
doing it effectively reduces the problem to one of consensus for
double-integrators. The latter problem has been researched
extensively in the literature (e.g.,
\cite{Ren3,Yu:2010ch,Tian}).

\subsection{Related work}

Typical coordination problems include attitude synchronization,
rendezvous, flocking, and formation control. For networks of single or
double-integrator systems, the rendezvous problem is referred to as
{\em consensus} or {\em agreement}, and it has been investigated by
many researchers, for
instance~\cite{Ren3,Yu:2010ch,Tian,beard,moreau,olfati,Spong,aut}.

A passivity-based solution of the attitude synchronization problem for
kinematic vehicle models is proposed in \cite{Hatanaka2}.
In~\cite{Ren1,Tayebi1,Abdessameud2}, the same problem is investigated
for dynamic vehicle models.  The proposed controllers do not require
measurements of the angular velocity, but they do require absolute
attitude measurements.  In~\cite{Nair}, the authors use the energy
shaping approach to design local and distributed controllers for
attitude synchronization. The same approach is adopted
in~\cite{Sarlette} to design two attitude synchronization controllers,
both local and distributed.  The first controller achieves
almost-global synchronization for directed connected graphs. However,
the controller design is based on distributed observers~\cite{scl},
and therefore requires auxiliary states to be communicated among
neighboring vehicles. It also employs an angular velocity dissipation
term that forces all vehicle angular velocities to zero in
steady-state. The second controller in~\cite{Sarlette} does not
restrict the final angular velocities, and does not require
communication, but it requires an undirected sensing graph, and
guarantees only local convergence.

The rendezvous problem for kinematic unicycles was solved
in~\cite{Lin} using time-varying feedbacks. The papers~\cite{Lin,
  SePaLe, Pappas,ConMorPraTos12} discuss the feasibility of achieving
various formations using local and distributed feedback for kinematic
unicycle models.  Dynamic unicycle models are considered in~\cite{Liu,
  Hawwary}. In~\cite{Liu}, a two-mode formation control is presented
in which the sensing graph has a spanning tree with a designated
leader vehicle as the root. Each vehicle, however, has access to the
acceleration of the leader through communication. The control strategy
requires a switch between two control modes designed to deal with
nonholonomic constraints in the system. The paper~\cite{Hawwary}
presents a local and distributed control law making dynamic and
kinematic unicycles converge to a common circle whose centre is
stationary and dependent on the initial configuration of the
vehicles. The spacing and ordering of unicycles on the circle is also
controlled. The problem is solved using a three step hierarchical
control based on a reduction theorem for the stabilization of sets.

The case of kinematic vehicles in three-space is investigated
in~\cite{Nair,Dong,Hatanaka,Scardovi}. The authors of~\cite{Nair,Dong}
consider the problem of full attitude and position synchronization,
but assume fully actuated vehicles. In~\cite{Scardovi}, the authors
propose distributed controllers to stabilize relative equilibria
which, as shown in \cite{Justh:2004in,Justh:2005un}, correspond to
parallel, circular or helical formations. Finally, in~\cite{Lee,
  Abdessameud} the authors consider formation control for dynamic,
underactuated vehicle models. However, the feedbacks are not local and
distributed. Also, in~\cite{Abdessameud} the sensing graph is assumed
to be undirected, and communication among vehicles is required, while
in~\cite{Lee} the graph is balanced, and it is assumed that each
vehicle has access to the thrust input of its neighbors, therefore
requiring once again communication between vehicles. Both approaches
in~\cite{Lee, Abdessameud} use a two-stage backstepping methodology in
which the first stage treats each vehicle as a point-mass system to
which a desired thrust is assigned. A desired thrust direction is then
extracted and backstepping is used to design a rotational control such
that vehicle rendezvous or formation control is achieved.  Our
previous work~\cite{RozMag14} investigates almost-global vehicle
rendezvous making use of a two-stage hierarchical methodology similar
to~\cite{Lee, Abdessameud}. In this approach, one can combine a
consensus controller for a network of double-integrators and an
attitude tracking controller satisfying certain assumptions to produce
a rendezvous controller for underactuated vehicles. However, this
approach requires that all vehicles can sense a common inertial vector
in their own body frame, which requires additional on-board sensors.
Moreover, the approach requires communication among vehicles.  The
solution presented in this paper overcomes all these limitations.  To
the best of our knowledge, a solution to the rendezvous control
problem for underactuated flying vehicles stated earlier has not yet
appeared in the literature.

\subsection{Organization of the paper.} 

We begin, in Section~\ref{sec:prelims}, by introducing some notation
and presenting basic notions of homogeneity of functions and stability
of sets. In Section~\ref{sec:modeling} we review the vehicle model. In
Section~\ref{sec:SCP} we formulate the rendezvous control problem. The
main result, Theorem~\ref{thm:2step}, is presented in
Section~\ref{sec:soln_pcp_2stage}, and its proof in
Section~\ref{sec:proof}.  In Section~\ref{sec:simulation}, we present
simulation results showing that the proposed solution is robust
against measurement errors, as well as force and torque
disturbances. Finally, in Section~\ref{sec:conclusion}, we end the
paper with some remarks. The proof of the main result relies on two
technical lemmas that are proved in the appendix.

\section{Preliminaries and notation}\label{sec:prelims}
We denote by $\Re_+$ the set of positive real numbers.  We use
interchangeably the notation $v=[v_1 \ \cdots \ v_n]^\top$ or
$(v_1,\ldots,v_n)$ for a column vector in $\Re^n$. We denote by $\one
\in \Re^n$ the vector $(1,\ldots,1)$.  If $v,w$ are vectors in
$\Re^3$, we denote by $v\cdot w:=v^\top w$ their Euclidean inner
product (also called the dot product), and by $\|v\|:=(v \cdot
v)^{1/2}$ the Euclidean norm of $v$.  If $v=(v_x,v_y,v_z)$, we define
\[
v^\times  := \left[ \begin{array}{rrr} 0 \hspace{1.6mm}  & \ \ -v_z & 
\ \ v_y \\ 
v_z & \ \ 0 \hspace{1.6mm} & \ \ -v_x \\ 
-v_y & \ \ v_x & \ \ 0 \hspace{1.6mm}
  \end{array}\right].
\]
One has that $v^\times w= v \times w$. Let
$\{e_1,e_2,e_3\}$ denote the natural basis of $\Re^3$ and $\SO(3):=\{M
\in \Re^{3 \times 3}: M^{-1} = M\trans, \det(M)=1\}$.  If $\Gamma$ is
a closed subset of a Riemannian manifold $\cX$, and $d: \cX \times \cX
\to [0,\infty)$ is a distance metric on $\cX$, we denote by
  $\|\chi\|_\Gamma:=\inf_{\psi \in \Gamma} d(\chi,\psi)$ the
  point-to-set distance of $\chi \in \cX$ to $\Gamma$. If
  $\varepsilon>0$, we let $B_\varepsilon(\Gamma) := \{\chi \in \cX :
  \|\chi\|_\Gamma < \varepsilon\}$ and by $\cN(\Gamma)$ we denote a
  neighborhood of $\Gamma$ in $\cX$. If $A,B \subset \cX$ are two
  sets, denote by $A \backslash B$ the set-theoretic difference of $A$
  and $B$. If $I = \{i_1,\ldots, i_n\}$ is an index set, the ordered
  list of elements $(x_{i_1},\ldots, x_{i_n})$ is denoted by $(x_j)_{j
    \in I}$.
  
Let $U, W$ be finite-dimensional vector spaces. A function $f: U \rightarrow W$ is {\it homogeneous} if, for all
$\lam>0$ and for all $x \in V$,
$f(\lam x)=\lam f(x)$.  A function $f: U \times V \rightarrow W$, $f(x,y)$ is {\it homogeneous with respect to $x$} if for all $\lam>0$ and for
all $(x,y) \in U \times V$, $f(\lam
x,y)=\lam f(x,y)$.

The following stability definitions are taken from~\cite{ElHMag12}.
Let $\Sigma: \dot \chi = f(\chi)$ be a smooth dynamical system with
state space a Riemannian manifold $\cX$. Let $\phi(t,\chi_0)$ denote
its local phase flow. Let $\Gamma \subset \cX$ be a closed set that is
positively invariant for $\Sigma$, i.e., for all $\chi_0 \in \Gamma$,
$\phi(t,\chi_0) \in \Gamma$ for all $t>0$ for which $\phi(t,\chi_0)$
is defined.
\begin{definition}
The set $\Gamma$ is {\it stable} for $\Sigma$ if for any $\varepsilon>0$,
there exists a neighborhood $\cN(\Gamma) \subset \cX$ such that, for
all $\chi_0 \in \cN(\Gamma)$, $\phi(t,\chi_0) \in B_\varepsilon(\Gamma)$,
for all $t >0$ for which $\phi(t,\chi_0)$ is defined. The set $\Gamma$
is {\it attractive} for $\Sigma$ if there exists neighborhood
$\cN(\Gamma) \subset \cX$ such that for all $\chi_0 \in \cN(\Gamma)$,
$\lim_{t\to\infty} \|\phi(t,\chi_0)\|_\Gamma =0$. The {\it domain of
  attraction of $\Gamma$} is the set $\{\chi_0 \in \cX:
\lim_{t\to\infty} \| \phi(t,\chi_0)\|_\Gamma =0\}$. The set $\Gamma$
is {\it globally attractive} for $\Sigma$ if it is attractive with
domain of attraction $\cX$. The set $\Gamma$ is {\it locally
  asymptotically stable (LAS)} for $\Sigma$ if it is stable and
attractive. The set $\Gamma$ is {\it globally asymptotically stable}
for $\Sigma$ if it is stable and globally attractive.
\hfill $\triangle$\end{definition}

Now consider a dynamical system $\Sigma(k): \dot \chi=f(\chi,k)$, in
which $k \in \Re^p$ is a vector of constant parameters (typically, control gains) and $f$ is a
smooth vector field with state space a Riemannian manifold.
\begin{definition} \label{def:GPS}
The set $\Gamma$ is {\it globally practically stable} for $\Sigma(k)$
if for any $\varepsilon > 0$, there exists a gain  $k^\star$ such
that $B_\varepsilon(\Gamma)$ has a subset which is globally asymptotically
stable for $\Sigma(k^\star)$.\hfill $\triangle$
\end{definition}

\section{Modeling}\label{sec:modeling}
We now return to the $i$-th robot depicted in
Figure~\ref{fig:vehicle}, with the aim of deriving its equations of
motion.  We fix a right-handed orthonormal inertial frame $\cI$,
common to all robots, and attach at the centre of mass of robot $i$ a
right-handed orthonormal body frame
$\cB_{i}=\{b_{ix},b_{iy},b_{iz}\}$, as depicted in the figure. We
denote by $(x_i,v_i)$ the inertial position and velocity of robot
$i$. We let $g$ denote the gravity vector in frame $\cI$.

We let $R_i$ be the $3 \times 3$ matrix whose columns are
the coordinate representations of $b_{ix}, b_{iy}, b_{iz}$ (in this
order) in frame $\cI$, so that $R_i \in \SO(3)$. The unit vector $q_i:=-R_i e_3$, depicted in
Figure~\ref{fig:vehicle}, is referred to as the {\em thrust direction vector} of
robot $i$, and the matrix $R_i$ is referred to as the {\em attitude}
of the robot. We assume that a thrust force $u_iq_i$ is applied at the centre of mass
of robot $i$. Notice that $u_i q_i$ has magnitude $u_i$, is directed
opposite to $b_{iz}$, and has constant direction in body frame
$\cB_i$.

Robot $i$ is assumed to have an actuation mechanism that
induces control torques $\tau_{ix}, \tau_{iy}, \tau_{iz}$ about its
body axes. We let $\tau_i := (\tau_{ix}, \tau_{iy}, \tau_{iz})$ be the
torque vector, and $\omega_i$ denote the angular velocity of the robot
with respect to frame $\cI$ (the unique vector in $\Re^3$ such that
$\dot R_i (R_i)^{-1} = \omega_i^\times$).

In this paper we adopt the convention that if $r \in \Re^3$ is an
inertial vector, the coordinate representation of $r$ in frame
$\cB_{i}$ is denoted by $r^i$, that is, $r^i:=R_i^{-1} r$. In
particular, the angular velocity of robot $i$ in its own body frame is
denoted by $\omega_i^i$. Finally, we use boldface symbols
  to denote reference quantities. For instance, $\wdes_i$ is
  the reference force for vehicle $i$ as in~\eqref{eq:dic} and $\OO_i$
  is the reference angular velocity for vehicle $i$ as
  in~\eqref{eq:tau}. The notation is summarized in
Table~\ref{table:abs_rel}.

Picking $(x_i,v_i,R_i,\omega_i^i)$ as state for robot $i$, we obtain
the equations of motion
\begin{equation}\label{eq:vehicle:translational}
\begin{aligned}
\dot{x}_{i} &= v_i, \\
m_i\dot{v}_i &= -u_iR_ie_3+m_ig, 
\end{aligned}
\end{equation}
\begin{equation}\label{eq:vehicle:rotational:so3}
\begin{aligned}
\dot R_i &= R_i \,(\omega_i^i)^{\times}, \\
J_i\dot{\omega}_i^i  &= \tau_i - \omega_i^i \times J_i\omega_i^i.
\end{aligned}
\end{equation}
In the above, $m_i$ is the mass of
robot $i$ and $J_i=J_i^\top$ is
its inertia matrix. We define the (inertial) relative
positions and velocities as $x_{ij}:=x_j - x_i$, $v_{ij}:=v_j -
v_i$. This model is standard
  and is widely used in the literature to model flying vehicles such
  as quadrotor helicopters. See, for instance,~\cite{Hua}. Sometimes
  researchers use alternative attitude representations, prominently
  quaternions~\cite{Abdessameud} or Euler
  angles~\cite{Mokhtari,Castillo}.  The
model~\eqref{eq:vehicle:translational}-\eqref{eq:vehicle:rotational:so3}
ignores aerodynamic effects such as drag and wind disturbances (such effects are included in~\cite{Hua}). It also ignores the
dynamics of the actuators.
\begin{table}[t]
\caption{Table of Notation}
\centering % used for centering table
\begin{tabular}{l l} 
\hline\hline %inserts double horizontal lines
Quantity & Description \\ [0.5ex] % inserts table
%heading
\hline % inserts single horizontal line
$m_i, \ J_i$ & mass and inertia matrix of robot $i$\\
$x_i \in \Re^3$ & inertial position of robot $i$\\ 
$v_i \in \Re^3$ & linear velocity of robot $i$\\
$R_i \in \SO(3)$ &  attitude of robot $i$\\
$\omega_i \in \Re^3$ & angular velocity of robot $i$\\
$q_i=-R_ie_3$ & thrust direction vector of robot $i$\\
$r^i = R_i^{-1} r$ & coord. repr. of  $r$  in frame $\cB_{i}$\\
$x_{ij}=x_j-x_i$ & rel. displacement of robot $j$ wrt robot $i$\\ 
$v_{ij}=v_j-v_i$ & rel. velocity of robot $j$ wrt robot $i$\\
%$R_{j}^i=R_i^{-1} R_j$ & rel. attitude of robot $j$ wrt robot $i$\\
%$q_{ij}=q_j-q_i$ & rel. thrust axis of robot $j$ wrt robot $i$\\
$\wdes_i \in \Re^3$ & reference force of robot $i$\\
$\OO_i \in \Re^3$ & reference angular velocity of robot $i$\\
${\cal N}_i$ & set of neighbors of robot $i$\\
$y_i=(x_{ij},v_{ij})_{j \in {\cal N}_i}$ & vector of rel. pos. and vel. available to robot $i$\\
%
%\hline %inserts single line
%\end{tabular}
%\centering % used for centering table
%\begin{tabular}{l l} 
%\hline\hline %inserts double horizontal lines
%Quantity & Description \\ [0.5ex] % inserts table
%%heading
%\hline % inserts single horizontal line
%$r^i = R_i^{-1} r$ & coord. repr. of  $r$  in frame $\cB_{i}$\\
%%
%$x_{ij}=x_j-x_i$ & rel. displacement of robot $j$ wrt robot $i$\\ 
%$v_{ij}=v_j-v_i$ & rel. velocity of robot $j$ wrt robot $i$\\
%%$R_{j}^i=R_i^{-1} R_j$ & rel. attitude of robot $j$ wrt robot $i$\\
%%$q_{ij}=q_j-q_i$ & rel. thrust axis of robot $j$ wrt robot $i$\\
%$\w_i= -u_i R_i e_3$ & applied thrust of robot $i$\\
%${\cal N}_i$ & set of neighbors of robot $i$\\
%$y_i=(x_{ij},v_{ij})_{j \in {\cal N}_i}$ & \\
\hline %inserts single line
\end{tabular}
\label{table:abs_rel} % is used to refer this table in the tex
\end{table}

\section{Rendezvous Control Problem}\label{sec:SCP}
We begin by defining the {\em sensor digraph}
$\mathcal{G}=(\mathcal{V},\mathcal{E})$, where $\mathcal{V}$ is a set
of nodes labelled as $\{1,\dots,n\}$, each representing a robot, and
$\mathcal{E}$ is the set of edges. An edge from node $i$ to node $j$
indicates that robot $i$ can sense robot $j$ (${\cal G}$ has no
self-loops).  A node is {\it globally reachable} if there exists a
path from any other node to it \footnote{For a graph $\mathcal{G}$, existence of a globally reachable node is equivalent to having a directed spanning tree in the reverse graph.}.

We denote by $\cN_i \subset {\cal V}$ the set of vehicles that robot
$i$ can sense. In a realistic scenario, $\cN_i$ is the set of robots
within the field of view of robot $i$. For instance, if each robot
mounted an omnidirectional camera, then one could define $\cN_i$ to
be the collection of robots that are within a given distance from
robot $i$. With such a definition, the sensor digraph ${\cal G}$ would
be state-dependent, making the stability analysis too hard at
present\footnote{Relatively little research has been done on
  distributed coordination problems with state-dependent sensor
  graphs. In this context, in the simplest case when the robots are
  modelled as kinematic integrators, it has been shown
  in~\cite{LinFraMag07} that the circumcentre law of Ando et
  al.~\cite{AndOasSuzYam99} preserves connectivity of the sensor graph
  and leads to rendezvous if the sensor graph is initially
  connected. Despite the simplicity of the robot model, the stability
  analysis in~\cite{LinFraMag07} is hard, and the control law is
  continuous but not Lipschitz continuous.}.

In light of the above, in this paper we assume that $\cN_i$ is
constant for each $i \in \{1,\ldots, n\}$ (and hence ${\cal G}$ is
constant as well).  If $j \in {\cal N}_i$, then we say that robot $j$
is a {\em neighbour} of robot $i$. If this is the case, then robot $i$
can sense the relative displacement and velocity of robot $j$ in its
own body frame, i.e., the quantities $x_{ij}^i, v_{ij}^i$.  Define the
vector $y_i:=(x_{ij}, v_{ij})_{j \in {\cal N}_i}$.  The relative
displacements and velocities available to robot $i$ are contained in
the vector $y_i^i:=(x^{i}_{ij}, v^{i}_{ij})_{j \in {\cal N}_i}$.  We
also assume that robot $i$ can sense its own angular velocity in its
own frame ${\cal B}_i$. To summarize, we have the definition below.

\begin{definition}\label{def:ldf}
A {\it local and distributed feedback} $(u_i,\tau_i)$ for robot $i$ is
a locally Lipschitz function of $y_i^i$ and $\omega_i^i$.
\hfill $\triangle$\end{definition}

The adjective {\it local} indicates that all quantities are
represented in the body frame of robot $i$, while {\it distributed}
indicates that only relative quantities with respect to neighboring
robots are accessible. In applications, a local and distributed feedback
for robot $i$ can be computed with on-board cameras and rate
gyroscopes.

We are now ready to define the Rendezvous Control Problem. 

\noindent {\em Rendezvous Control Problem:} Consider
system~\eqref{eq:vehicle:translational},~\eqref{eq:vehicle:rotational:so3}, and define the {\em rendezvous manifold}
\begin{equation} \label{eq:rendez_set}
\begin{aligned}
\Gamma := &\left\{(x_i,v_i,R_i,\omega_i^i)_{i \in \{1, \dots, n\}} \in \Re^{3n} \times \Re^{3n} \times \SO(3)^n \times \Re^{3n} \right.\\
&: \left. x_{ij}=v_{ij}=0, \ \forall i,j \right\}.
\end{aligned}
\end{equation}
Find, if possible, local and distributed feedbacks $(u_i,\tau_i)_{i
  =1,\ldots,n}$ that globally practically stabilize $\Gamma$. \hfill $\triangle$

The goal of the rendezvous control problem is to achieve
synchronization of the robot positions and velocities to any desired
degree of accuracy from any initial configuration.

\section{Solution of the Rendezvous Control Problem}\label{sec:soln_pcp_2stage}

\begin{definition}\label{defn:CC}
Consider a collection of $n$ double-integrators
\begin{equation} \label{eq:CC_sys}
\begin{aligned}
\dot x_i&=v_i\\
\dot v_i&=\w_i, \ i=1 \dots n,
\end{aligned}
\end{equation}
where $f_i$ is the control input of subsystem $i$.
Suppose the double-integrators have the
same sensor digraph ${\cal G}$ as the underactuated robots of
Section~\ref{sec:modeling}. A feedback $\wdes_i(y_i)$, $i=1 \dots n$,
is a \textit{double-integrator consensus controller} if $\wdes_i$ has
the form
\begin{equation}\label{eq:dic}
\wdes_i(y_i)=\sum_{j \in {\cal N}_i} \Big( a_{ij}x_{ij}+ b_{ij}v_{ij}
\Big), \ i=1,\ldots, n,
\end{equation}
with $a_{ij}, \ b_{ij} \in \Re$ and if, setting $\w_i=\wdes_i(y_i)$
in~\eqref{eq:CC_sys}, the set
\[
\left\{(x_i,v_i)_{i\in \{1 \dots n\}} \in \Re^{3n} \times \Re^{3n}:
x_{ij}=0, v_{ij}=0, \ \forall i,j \right\}
\]
is globally asymptotically stable for~\eqref{eq:CC_sys}. 
\hfill $\triangle$\end{definition}

Ren et al. in~\cite[Theorems 4.1, 4.2]{Ren3} and Yu et
al. in~\cite[Theorem 1]{Yu:2010ch} have shown that a double-integrator
consensus controller exists if and only if the sensor digraph ${\cal
  G}$ has a globally reachable node.  Now the main result of this
paper.

\begin{theorem}\label{thm:2step}
If the sensor digraph ${\cal G}$ has a globally reachable node, then
the rendezvous control problem is solvable for
system~\eqref{eq:vehicle:translational}-\eqref{eq:vehicle:rotational:so3},
and a solution is given as follows.  Let $\wdes_i$, $i=1,\ldots,n$, be
a double-integrator consensus controller. The local and distributed
feedback,
%\begin{equation} \label{eq:CCP}
\begin{align}\label{eq:CCP}
\begin{split}
u_i =& - m_i\wdes_i(y_i^i)\cdot e_3, \\ 
\tau_i=&\omega_i^i \times J_i
  \omega_i^i - \kq J_i\left((\omega_i^i \times \wdes_i(y_i^i))\times
  e_3\right) \end{split}\\
  &-\kq^2 \ktilde \left[\omega_i^i -
    \kq(\wdes_i(y_i^i)\times e_3)\right], \ i=1 \dots n, \nonumber
\end{align}
%\end{equation}
where $\kq, \ktilde>0$ are control parameters, makes the rendezvous
manifold~\eqref{eq:rendez_set} globally practically stable. In
particular, for any $\varepsilon>0$, there exist
$\kq^\star,\ktilde^\star>0$ such that for all $\kq>\kq^\star$,
$\ktilde > \ktilde^\star$, the set $B_\varepsilon(\Gamma)$ has a globally
asymptotically stable subset.
\end{theorem}

The proof of Theorem~\ref{thm:2step} is presented in
Section~\ref{sec:proof}.  

\subsection*{Explanation of  proposed controller}
Returning to the block diagram of Figure~\ref{fig:block_diagram}, we
now explain in detail the operation of its two nested loops. We begin with the observation that a
double-integrator consensus controller $\wdes_i(y_i)$, $i=1 \dots n$,
for system~\eqref{eq:CC_sys} also makes the systems
\begin{equation} \label{eq:dcg}
\begin{aligned}
\dot x_i&=v_i\\
\dot v_i&=\w_i+g
\end{aligned}
\end{equation}
rendezvous, since the addition of the gravity vector $g$ does not affect the relative dynamics.
Now compare system~\eqref{eq:dcg} to the translational dynamics of the
flying robots,
\begin{equation}\label{eq:translational}
\begin{aligned}
\dot x_i &=v_i \\
\dot v_i &=-\frac{1}{m_i}u_i R_i e_3 + g.
\end{aligned}
\end{equation}
If it were the case that $\w_i = -(1/m_i) u_i
  R_i e_3$, systems~\eqref{eq:dcg} and~\eqref{eq:translational} would
  be identical.  Then, setting $-u_i R_i e_3=
    m_i \wdes_i$ in~\eqref{eq:translational} would solve the
rendezvous problem. Inspired by this observation, the outer loop of
the block diagram in Figure~\ref{fig:block_diagram} assumes that $-u_i R_i e_3$ is the control input
of~\eqref{eq:translational} and computes a desired double-integrator
force $m_i \wdes_i$ which becomes a reference signal for the inner
loop.

\begin{figure}[t]
\centerline{\includegraphics[width=0.25\textwidth]{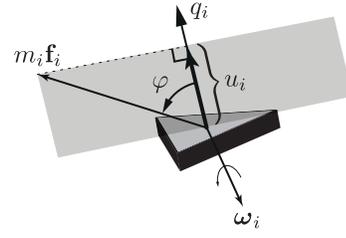}}
\caption{Illustration of the control input $u_i$ and reference
  angular velocity $\OO_i$ in~\eqref{eq:CCP}.}
\label{fig:vehicle_detailed}
\end{figure}

We now explore in more detail the operation of the inner loop. First we observe that since  $\wdes_i$ is a linear function, we have 
$R_i \wdes_i(y_i^i) = \wdes_i(R_i y_i^i) = \wdes_i(y_i)$. Moreover, using the fact that 
dot products are invariant under rotations, we have
\[
u_i = -m_i\wdes_i(y_i^i) \cdot e_3 = m_i(R_i\wdes_i(y_i^i)) \cdot (-R_i e_3) = m_i\wdes_i(y_i) \cdot q_i,
\]
where $q_i$ is the thrust direction vector. Thus, the thrust magnitude is the
projection of the desired thrust $m_i\wdes_i$ onto the thrust direction vector---see
Figure~\ref{fig:vehicle_detailed}. Now let
$\OO^i_i(y_i^i)=\kq\left(\wdes_i(y^i_i) \times e_{3} \right)$. Then we
have
\begin{equation}\label{eq:tau}
\begin{aligned}
\tau_i \hspace{-0.5mm}=&\hspace{-0.5mm} \omega_i^i \times J_i \omega_i^i - \kq J_i\left((\omega_i^i
\times \wdes_i(y_i^i))\times e_3\right)\\
&-\kq^2 \ktilde \left(\omega_i^i - \OO_i^i(y_i^i)\right).
\end{aligned}
\end{equation}
We will show in the proof of Theorem~\ref{thm:2step} that the
  torque inputs $\tau_i$ make $\omega_i^i$ converge to an arbitrarily
  small neighborhood of $\OO_i^i$, $i=1,\ldots, n$. Thus, $\OO_i^i$
can be seen as a reference angular velocity for the inner loop. Using the fact
that, for all $a, b \in \Re^3$ and all $R \in \SO(3)$, $R(a \times b)
= (Ra) \times (Rb)$, we have
\[
\begin{aligned}
\OO_i &= R_i \OO_i^i = R_i \kq\left(\wdes_i(y^i_i) \times e_{3} \right) = \kq
\left( (R_i\wdes_i(y^i_i)) \times (R_i e_3) \right)\\ 
&= \kq
(\wdes_i(y_i) \times -q_i) 
=\kq (q_i \times \wdes_i(y_i)).
\end{aligned}
\]
Thus $\OO_i$ is perpendicular to the plane formed by the thrust direction vector $q_i$ and the desired thrust
force $m_i\wdes_i$---see Figure~\ref{fig:vehicle_detailed}. Since the
angular velocity vector identifies an instantaneous axis of rotation,
it follows that if $\omega_i= \OO_i$, then robot $i$ rotates about
$\OO_i$ according to the right-hand rule. Referring to
Figure~\ref{fig:vehicle_detailed}, we see that such a rotation closes
the gap between $u_iq_i$ and
$m_i\wdes_i$, and the speed of rotation is
proportional to $\sin \varphi$, where $\varphi$ is the angle between
$u_iq_i$ and $m_i\wdes_i$ marked in the figure. When the gap is
closed, we have $u_i = \|m_i\wdes_i\|$, $q_i = m_i\wdes_i / \|m_i\wdes_i\|$, and
thus $u_i q_i =m_i\wdes_i$.
In conclusion, the inner loop assigns $(u_i,\tau_i)$ to make
$\omega_i$ approximately converge to $\OO_i$, so that
$u_iq_i= -u_i R_i e_3$ approximately
converges to $m_i\wdes_i$, which is computed by the outer loop.

While the intuition behind the proposed controller is simple, the
proof that the interplay between the two nested loop results in global
practical stability of the rendezvous manifold is rather delicate, and
it crucially relies on the homogeneity of the functions $\wdes_i$,
$i=1,\ldots, n$.
\begin{remark}
Theorem~\ref{thm:2step} proves global practical stability of the rendezvous manifold $\Gamma$. The reason that the stability is practical and not asymptotic is roughly as follows. In order to achieve rendezvous of the rigid bodies, $u_iq_i$ is driven approximately to $m_i\wdes_i$. What's important is not so much the difference in magnitude of these vectors  but rather the difference in angle between them. In Figure 3, one can see that $\OO_i$ acts to reduce this angle with a rate proportional to the magnitude of $\OO_i$. Since $\OO_i$ is a linear function of $\wdes_i$,  as the robots approach consensus $\OO_i$ converges to zero at the same rate as  $\wdes_i$. This leads to increasing inaccuracy in closing the gap between the vectors $u_iq_i$ and $m_i\wdes_i$ insomuch that in a very small neighborhood of rendezvous, $\OO_i$ is so small that it fails to make the translational dynamics act as double integrators. More detailed reasoning is provided in Remark~\ref{rem:gas}.
\end{remark}

\subsection*{Features of the proposed controller}

(i) \hspace{1pt} The proposed controller has a number of
advantages over our previous work
in~\cite{RozMag14}. Unlike~\cite{RozMag14}, the inner control loop
does not require any derivatives of the reference thrust force
$\wdes_i$. In~\cite{RozMag14}, the large expressions resulting from
such derivatives pose difficulty in real-time computation of the
control law. More importantly, the computation of such derivatives
requires communication between neighboring robots, a problem that
has been overcome in the present approach.  The approach
in~\cite{RozMag14} requires that robots have access to a common
inertial vector. This requirement is absent in this paper.

(ii) \hspace{1pt} The feedback of Theorem~\ref{thm:2step} is
static. It does not depend on dynamic compensators that require
communication between neighboring robots.

(iii) \hspace{1pt} The feedback of Theorem~\ref{thm:2step} is
local and distributed in the sense of
Definition~\ref{def:ldf}. Interestingly, it does not require sensing
of relative attitudes, which can be computed using on-board cameras,
but are harder to compute than relative displacements.

(iv) \hspace{1pt} On the rendezvous manifold $\Gamma$ there is
no prespecified thrust direction $q_i$ for robot $i$ and the robot
thrust directions do not need to align at rendezvous. This is
desirable if one wants to employ the proposed controller in a
hierarchical control setting to enforce additional control
specifications.

\section{Proof of Theorem \ref{thm:2step}}\label{sec:proof}
The feedback in~\eqref{eq:CCP} is local and distributed because it is
a smooth function of $y_i^i$ and $\omega_i^i$ only. By Theorems 4.1
and 4.2 in~\cite{Ren3} (or Theorem 1 in~\cite{Yu:2010ch}), if ${\cal
  G}$ has a globally reachable node then there exists a
double-integrator consensus controller, and the
feedback~\eqref{eq:CCP} is well-defined. We need to show that it
renders the rendezvous manifold $\Gamma$ in~\eqref{eq:rendez_set}
globally practically stable. We begin by expressing the translational
portion of the dynamics in coordinates relative to robot $1$, i.e., in
terms of the variables $(x_{1j},v_{1j})_{j=2,\ldots,n}$,
\begin{align}
\begin{split} 
	\dot x_{1j} &= v_{1j}, \label{eq:rel_coords_trans}\\
	\dot v_{1j} &= -\frac{1}{m_j} R_j e_3 u_j +\frac{1}{m_1}R_1 e_3u_1, \quad j=2,\ldots, n, \\ 
\end{split}
\\
\begin{split} 
&	\dot R_{i} = R_i(\omega_i^i)^\times \label{eq:rel_coords_rot},\\
 & J_i\dot \omega_i^i = \tau_i-\omega_i^i \times J_i\omega_i^i, \quad i=1,\ldots,n.
\end{split}
\end{align}
Since all relative states $(x_{ij},v_{ij})$ can be expressed in terms
of the variables above through the identity $(x_{ij},v_{ij}) =
(x_{1j}-x_{1i},v_{1j}-v_{1i})$, perfect rendezvous occurs if and only
if the vector $(x_{1j},v_{1j})_{j=2,\ldots,n}$ is zero. Denoting, 
\[
\begin{aligned}
&\Chi:=(x_{1j},v_{1j})_{j=2,\ldots,n} \in \sX:=\Re^{3 (n-1)} \times
\Re^{3 (n-1)}, \\
&R:=(R_1,\ldots, R_n) \in \sR:=\SO(3)^n, \\
&\omega:=(\omega_1^1,\ldots, \omega_n^n) \in \sO :=\Re^{3n},
\end{aligned}
\]
the new collective state is $(\Chi,R,\omega) \in \sX \times \sR \times
\sO$.  The meaning of the new state is this: $\Chi$ contains all
translational states (positions and velocities) relative to robot $1$,
$R$ contains all the attitudes, and $\omega$ contains all body frame
angular velocities.  The rendezvous manifold in new coordinates is the
set $\{(\Chi,R,\omega) \in \sX \times \sR \times \sO: \Chi =0\}$.

Due to the identity $(x_{ij},v_{ij}) = (x_{1j}-x_{1i},v_{1j}-v_{1i})$,
the vector $y_i=(x_{ij},v_{ij})_{j \in \cN_i}$ is a linear function of
$\Chi$ which we will denote $y_i=h_i(X)$. Similarly, the vector
$y_i^i = (x_{ij}^i,v_{ij}^i)_{j \in \cN_i}$ is a function of $\Chi$
and $R$, linear with respect to $\Chi$. We will denote this function
$y_i^i = h_i^i(\Chi,R)$.

Using the definitions above, we may now express $\wdes_i^i(y_i^i)$ and
$\OO_i(y_i^i)=\kq (\wdes_i(y_i^i) \times e_3)$ (the latter function
was discussed in Section~\ref{sec:soln_pcp_2stage}) in terms of
states.  Accordingly, we define $\gdes_i: \sX \to \Re^3$, $\gdes_i^i:
\sX \times \sR \to \Re^3$ and $\OO : \sX \times \sR \to \sO$ as
follows:
\begin{equation}\label{eq:g_and_OO}
\begin{aligned}
\gdes_i(\Chi)&:=\wdes_i \circ h_i(\Chi), \\
\gdes_i^i(\Chi,R)&:=R_i^{-1} \gdes_i(\Chi) = \wdes_i \circ
h_i^i(\Chi,R), \\
\OO(\Chi,R)&:=\big(
\OO_i(h_i^i(\Chi,R)) \big)_{i=1,\ldots,n}.
\end{aligned}
\end{equation}
We remark that $\gdes_i$ is linear and $\gdes_i^i$ is linear with
respect to its first argument. The second identity in the definition
of $\gdes_i^i$ is due to the linearity of $\wdes_i$.

Finally, we define the rendezvous manifold in new coordinates,
\begin{equation}\label{eq:Gamma_star}
\Gamma^\star:=\{(\Chi,R,\omega) \in \sX \times \sR \times \sO: \Chi=0 \}.
\end{equation}
We will prove that $\Gamma^\star$ is globally practically stable,
which will imply that $\Gamma$ is globally practically stable as well.

\subsection{Lyapunov function}

Consider the $n$ double-integrators~\eqref{eq:CC_sys} with control $\wdes_i$
in~\eqref{eq:dic}, expressed in $\Chi$ coordinates:
\begin{equation}\label{eq:consensus_closed_loop}
\begin{aligned}
\dot x_{1j} &= v_{1j} \\
\dot v_{1j} &= \wdes_j(y_j) - \wdes_1(y_1) = \gdes_j(\Chi) -
\gdes_1(\Chi), \quad j=2,\ldots, n.
\end{aligned}
\end{equation}
By Definition~\ref{defn:CC}, the origin of this linear time-invariant
system is globally asymptotically stable. Thus, there exists a
quadratic Lyapunov function $V : \sX \to \Re$, $V(\Chi) =
\Chi^\top P \Chi$, where $P$ is a symmetric positive definite
matrix, such that the derivative of $V$ along the vector field
in~\eqref{eq:consensus_closed_loop} is negative definite.

Let $\Je \in \Re^{3n \times 3n}$ be the block-diagonal matrix with the $i$-th block equal to $J_i$, and consider the function $W:\sX \times \sR
\times \sO \rightarrow \Re$ defined as
\begin{equation}\label{eq:W}
 W(\Chi,R,\omega)=\alpha W_{\tran}(\Chi)+ W_{\rot}(\Chi,R,\omega),
\end{equation}
where $\alpha>0$ is a parameter to be assigned
  later and
\[
\begin{aligned}
W_{\tran}(\Chi) &=  \sqrt{V(\Chi)}+\frac{1}{2} V(\Chi), \\
W_{\rot}(\Chi,R,\omega) =&\sum_{i=1}^n \gdes_i^i (\Chi,R) \cdot
 e_3\\
 &+\frac{1}{2} (\omega - \OO(\Chi,R))\trans \Je (\omega - \OO(\Chi,R)).
\end{aligned}
\]

\begin{lemma}\label{lem:saturation}
Consider the continuous function $W$ defined in~\eqref{eq:W}.
Then 
\[
\alpha^\star:=\sup_{(X,R) \in \sX\backslash\{0\} \times
    \sR} \sum_i |\gdes_i^i(X/\sqrt{V(X)},R) \cdot e_3| < \infty,
\]
and for all $\alpha > \alpha^\star$, the following properties hold:
\begin{enumerate}[(i)]
\item $W \geq 0$ and $W^{-1}(0) \subset \Gamma^\star$.
\item For all $c>0$, the sublevel set $W_c :=\{(\Chi,R,\omega):
  W(\Chi,R,\omega) \leq c\}$ is compact.
\item For all $\varepsilon>0$, there exists
$\delta>0$ such that $W_\delta \subset B_\varepsilon(\Gamma^\star)$.
\end{enumerate}
\end{lemma}
The proof is in the appendix. 

From now on we assume $\alpha > \alpha^\star$. In light of the lemma,
if we show that $W$ is nonincreasing outside a certain compact region
of the state space, then all trajectories
of~\eqref{eq:rel_coords_trans}-\eqref{eq:rel_coords_rot} with
feedback~\eqref{eq:CCP} are bounded, ruling out finite escape
times. Moreover, in light of part (iii) of the lemma, to prove that
$\Gamma^\star$ is practically stable it suffices to prove that for
every $\delta>0$, there exists a gain vector $(\kq,\ktilde)$ such that
$W_\delta$ is globally asymptotically stable. For this, we need to
show that $W \geq \delta \implies \dot W<0$.

\subsection{Coordinate transformation}

\begin{figure}[t]
\centerline{\includegraphics[width=0.25\textwidth]{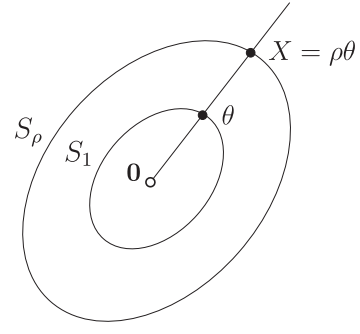}}
\caption{Illustration of the sets $S_1$ and $S_{\lam}$.}
\label{fig:proof}
\end{figure}

We now construct a coordinate transformation on the translational
states $\Chi$ that leverages the homogeneity property of
$\wdes_i$. Return to the Lyapunov function $V(\Chi) = \Chi\trans P
\Chi$ associated with the double-integrator consensus
controller. Since $V$ is a positive definite quadratic form, its level
sets are compact and convex. Consider the level set $S_1 := \{\Chi \in \sX: V(\Chi) =\Chi^\top P \Chi = 1 \},$
and for $\lam>0$, let $S_{\lam}$ denote the set $S_{\lam}:=\{\Chi \in \sX: \Chi = \lam \prj, \prj \in S_1 \}$. The sets $S_1$ and $S_{\lam}$ are depicted in
Figure~\ref{fig:proof}.  
By convexity of $S_1$, any point $\Chi \in
\sX \backslash \{0\}$, can be uniquely represented as $
\Chi = \lam \prj, \quad \lam \in \Re_+, \ \prj \in S_1$,
where $\lam = \sqrt{\Chi^\top P \Chi}$ and $\prj = \Chi /
\lam$. In the above decomposition, one can think of $\lam$ as a
scaling factor determining the size of the neighborhood of zero where
$\Chi$ belongs, while $\prj$ is a shape variable determining the
relative positions and velocities of the robots modulo scaling.
We use this construction to transform the coordinates of the relative
translational states in $\Chi$ as follows. Define the map 
\[
\begin{aligned}
& F: \sX \backslash \{0\} \times \sR \times \sO \to \Re_+ \times S_1
  \times \sR \times \sO, \\
& F(\Chi,R,\omega)=(\lam,\prj,R,\omega), \ \lam:=\sqrt{V(\Chi)},
  \prj:= \Chi / \sqrt{V(\Chi)}.
\end{aligned}
\]
Clearly $F$ is a smooth bijection. Moreover its inverse $
F^{-1}(\lam,\prj,R,\omega) = (\lam \prj,R,\omega)$ is smooth as well, so $F$ is a diffeomorphism\footnote{$F$ is a
  diffeomorphism of smooth manifolds. The set $S_1$ is diffeomorphic
  to the unit sphere of dimension $6(n-1)-1$. All other sets involved
  in the Cartesian products are smooth manifolds}.  The new state is
$(\lam,\prj,R,\omega) \in \Re_+ \times S_1 \times \sR \times \sO$.
Rendezvous in these coordinates would correspond to having
$\lam=0$, which is outside of the image of $F$. This is not a
problem though, since we want to show practical stability of the
rendezvous manifold, for which it suffices to show that $\lam$ can
be made arbitrarily small.

Having defined a coordinate transformation, our next objective is to
represent the Lyapunov function candidate $W$ in new coordinates. The
new representation is $\hat W = W \circ F^{-1}$, which amounts to
simply replacing $\Chi$ by $\lam \prj$. Doing so we obtain
\[
\hat W(\lam,\prj,R,\omega) = \alpha \hat W_\tran(\lam)+ \hat
W_\rot(\lam,\prj,R,\omega),
\]
where $\hat W_\tran (\lam) =\lam +
\frac{\lam^2} 2$, 
\[
\begin{aligned}
\hat W_\rot (\lam,\prj,R,\omega) =&\lam \sum_{i=1}^n
\gdes_i^i(\prj,R) \cdot e_3 \\
&+ \frac 1 2 \big( \omega - 
\OO(\lam\prj,R)  \big)\trans \Je \big(\omega -
\OO(\lam\prj,R) \big).
\end{aligned}
\]
In writing the above, we used the identity $\lam = \sqrt{V(\chi)}$
and the fact that the function $\gdes_i^i(\Chi,R)$ is linear with
respect to $\Chi$, implying that $\gdes_i^i(\lam \prj,R) = \lam
\gdes_i^i(\prj,R)$. In what follows, we let $\hat W_\delta :=
\{(\lam,\prj,R,\omega) \in \Re_+ \times S_1 \times \sR \times \sO:
\hat W(\lam,\prj,R,\omega) < \delta\}$. Thus, $\hat W_\delta = F(W_\delta)$.

\subsection{Stability analysis}

Let $\delta>0$ be arbitrary. We have $\hat W \leq \alpha(\rho+\rho^2/2) + \lam
  \sup_{(\prj,R)} |\gdes_i^i(\prj,R) \cdot e_3| + (1/2) (\omega -
  \OO)\trans \Je (\omega - \OO)$. Using the definition of
  $\alpha^\star$ in Lemma~\ref{lem:saturation} and the fact that
  $\alpha > \alpha^\star$, we get 
\[
\hat W \leq \alpha(2\lam+\lam^2/2) + (1/2) (\omega - \OO)\trans \Je
(\omega - \OO).
\]
It readily follows that there exists $\varrho\in(0,1)$ such that
\[
\Lambda_{\varrho} :=\{(\lam,\prj,R,\omega) : \lam \in
(0,\varrho), \|\omega - \OO(\lam\prj,R)\|^2 < \varrho \}
\subset \hat W_\delta.
\]
%% \textcolor{red}{This seems incorrect since the inequality is $\geq$. We should use the inequality $\hat W
%% \leq 2\alpha \lam + (1/2) \alpha \lam^2 + (1/2)(\omega - \OO )\trans \Je (\omega -
%% \OO)$}
We will show that there exist $\alpha>0$ and a gain vector
$(\kq,\ktilde)$ such that $\dot {\hat W} <0$ outside the set
$\Lambda_{\varrho}$. This will imply that $\hat W \geq \delta \implies
\dot {\hat W} <0$, proving that $\hat W_\delta$ is globally
asymptotically stable.

\begin{lemma}\label{lem:Wt}
Consider the closed-loop
system~\eqref{eq:rel_coords_trans}-\eqref{eq:rel_coords_rot} with
feedback~\eqref{eq:CCP}. If $\kq > 1$, then there exist scalars
$M_1,\ldots, M_4>0$ such that the derivatives of $\lam$ and $\hat
W_\rot$ along the closed-loop system in $(\lam,\prj,R,\omega)$
coordinates satisfy the following inequalities:
\begin{align}\label{eqn:inequalities}
\begin{split}
 \dot\lam\leq&\lam\left[
-M_{2} + M_1 \sum_{i=1}^n \|\gdes_i^i(\prj,R) \times e_3\|\right], \\
 \dot {\hat W}_\rot \leq& \lam^2 \sum_{i=1}^n
\left[-\kq\|\gdes_i^i(\prj,R) \times e_3 \|^2 \right.  \end{split}\\
&\left. + \frac{M_4}{\ktilde}\right]+ \lam M_{3} - \frac{\kq^2\ktilde}{2}\|\omega- \OO(\lam\prj,R)\|^2.\nonumber
\end{align}
\end{lemma}
The proof  is in the appendix. 

From now on we let $\kq>1$. Using the inequalities in
Lemma~\ref{lem:Wt}, we get
\[
\begin{aligned}
\dot{\hat W} \leq& (\lam+ \lam^2)\left[-\alpha M_{2}+ \alpha M_1
  \sum_{i=1}^n \|\gdes_i^i(\prj,R) \times e_3\|\right] \\
&+\lam^2 \sum_{i=1}^n \left[-\kq\|\gdes_i^i(\prj,R) \times e_3 \|^2 +
  \frac{M_4}{\ktilde}\right]\\
&+ \lam M_{3}- \frac{\kq^2\ktilde}{2}\|\omega-\OO(\lam\prj,R)\|^2.
\end{aligned}
\] 
Denote $\bb_i(\prj,R):=\|\gdes_i^i(\prj,R) \times e_3\|$, and
$\bb(\prj,R):=(\bb_1(\prj,R),\ldots,\bb_n(\prj,R))$. For notational
convenience, we omit the arguments of the functions $\bb$ and
$\OO$. With these definitions, the inequality above may be rewritten
as
\[
\begin{aligned}
\dot{\hat W} \leq& (\lam+ \lam^2)\left( -\alpha M_{2}+ \alpha M_1
\one\trans \bb \right) \hspace{-0.5mm} +\hspace{-0.5mm} \lam^2 \hspace{-0.5mm} \left (\hspace{-0.5mm} -\kq\|\bb\|^2 \hspace{-1mm} +\hspace{-1mm}  \frac{M_4
  n}{\ktilde}\right) \\
&+ \lam M_{3}- \frac{\kq^2\ktilde}{2}\|\omega-\OO\|^2.
\end{aligned}
\] 
For every $\ktilde > n M_4 / M_3$, we have
\[
\begin{aligned}
\dot{\hat W}\leq &(\lam+ \lam^2)\left( -\alpha M_{2} + M_{3}+
\alpha M_1 \one\trans\bb \right) - \lam^2 \kq \|\bb\|^2 \\
&-
\frac{\kq^2\ktilde}{2}\|\omega-\OO\|^2.
\end{aligned}
\]
If we further pick $\alpha > \max\{\alpha^\star, 3 M_3 / M_2\}$, we
have
\[
\begin{aligned}
\dot{\hat W} \leq& (\lam+ \lam^2)\left( -2M_{3}+ \alpha M_1
\one\trans\bb\right) -\lam^2 \kq \|\bb\|^2 \\
&-
\frac{\kq^2\ktilde}{2}\|\omega-\OO\|^2.
\end{aligned}
\]

Splitting the term $-\lam^2 \kq \|\bb\|^2$ into two parts and
collecting terms for $\lam$ and $\lam^2$, we obtain
\[
\begin{aligned}
\dot{\hat W} \leq& 
\lam^2 \left(-2M_{3}+ 
\alpha M_1 \one\trans\bb- \frac{\kq}{2} \|\bb\|^2\right) \\
&+ \lam\left(-2M_{3}+ \alpha M_1 \one\trans \bb
-\lam\frac{\kq}{2} \|\bb\|^2\right) \\
&- 
\frac{\kq^2\ktilde}{2}\|\omega -\OO\|^2.
\end{aligned}
\]
Consider now the expression
\[
M_3 - \alpha M_1 \one\trans \bb + \frac {\kq \varrho} 2 \|
\bb\|^2  \hspace{-1mm}=\hspace{-1mm} \big[ \one\trans \ \ \bb\trans\big] 
\hspace{-1mm} 
\begin{bmatrix}
  \frac{M_3} n I  \hspace{-1mm} &\hspace{-1mm} -\alpha \frac{M_1}{2} I \\
-\alpha \frac{M_1}{2} I \hspace{-1mm} &\hspace{-1mm} \frac{\kq \varrho}2 I 
\end{bmatrix}
\hspace{-1mm} 
\begin{bmatrix} 
\one \\ \bb
\end{bmatrix}.
\]
If $\kq >2 n(\alpha M_1 /2)^2 /(\varrho M_3)$, the above
quadratic form is positive definite, implying that
\begin{equation}\label{eq:quadratic}
M_3 - \alpha M_1 \one\trans \bb + \frac {\kq \varrho} 2 \|
\bb\|^2 \geq 0.
\end{equation}
Since $\varrho<1$, we also have $M_3 - \alpha M_1 \one\trans \bb +
(\kq/2) \| \bb\|^2 \geq 0$. Using the latter inequality, we get a
further upper bound for $\dot {\hat W}$,
\begin{equation}\label{eq:dotW_bound}
\begin{aligned}
\dot {\hat W} \leq& -\lam^2M_{3}+\lam\left(-2M_{3}+ \alpha M_1
\one\trans\bb - \lam \frac{\kq}{2}\|\bb\|^2 \right)\\
&-
\frac{\kq^2\ktilde}{2}\|\omega-\OO\|^2.
\end{aligned}
\end{equation}
Using~\eqref{eq:dotW_bound}, we now prove that outside
$\Lambda_{\varrho}$, $\dot {\hat W}<0$. In other words,
when either $\lam \geq \varrho$ or $\|\omega-\OO\|^2
\geq \varrho$, $\dot {\hat W}<0$.

\begin{remark}\label{rem:gas}
If the derivative $\dot {\hat W}$ were negative definite, then the rendezvous manifold $\Gamma^\star$ would be globally asymptotically stable. However, this is not guaranteed in~\eqref{eq:dotW_bound}. The reason is as follows. Suppose $\lam$ is very small and $\|\omega-\OO\|=0$. Then all terms multiplied by $\rho^2$ become negligible and what remains in~\eqref{eq:dotW_bound} is, $\dot {\hat W} \leq \lam\left(-2M_{3}+ \alpha M_1\one\trans\bb \right)$. As we have no control over the value of the constants $M_1$ and $M_3$ in the equation above, $\dot {\hat W}$ can be greater than zero if the second term dominates the first.  
\end{remark}

Suppose first that $\lam \geq \varrho$. Then
from~\eqref{eq:dotW_bound} we have
\[
\begin{aligned}
\dot {\hat W} \leq& -\lam^2M_{3}+\lam\left(-2M_{3}+ \alpha M_1
\one\trans\bb - \frac{\kq \varrho}{2}\|\bb\|^2 \right)\\
&-
\frac{\kq^2\ktilde}{2}\|\omega-\OO\|^2.
\end{aligned}
\]
By inequality~\eqref{eq:quadratic} we conclude that
\[
\dot {\hat W} \leq -\lam^2 M_{3}-\lam M_{3} 
-\frac{\kq^2\ktilde}{2}\|\omega-\OO\|^2 <0.
\]

Next, suppose that $\|\omega-\OO\|^2 \geq \varrho$. Then
from~\eqref{eq:dotW_bound},
\[
\begin{aligned}
\dot {\hat W} &\leq -\lam^2M_{3} + \lam \alpha M_1 \one\trans\bb
- \frac{\kq^2\ktilde}{2} \varrho \\
&\leq -\lam^2M_{3} + \lam \alpha M_1 M_5
- \frac{\kq^2\ktilde}{2} \varrho,
\end{aligned}
\]
where $M_5:=\max_{(\prj,R)\in S_1 \times \sR} \{ \one\trans\bb(\prj,R)
\}$. The maximum exists because $\bb$ is continuous and $ S_1 \times
\sR$ is a compact set.  If $\ktilde > (\alpha M_1 M_5 /\kq)^2/\varrho$
then $\dot {\hat W}<0$.

We have therefore proved that, if $\alpha>
\max\{\alpha^\star,3M_3/M_2\}$, $\kq > \max\{1,2 n(\alpha M_1 /2)^2
/(\varrho M_3)\}$, and $\ktilde>\max\{n M_4 / M_3,(\alpha
M_1 M_5/\kq)^2/\varrho\}$, then $\hat W>\delta$
implies that $\dot {\hat W}<0$. Therefore, for any initial
condition, the solution
of~\eqref{eq:rel_coords_trans}-\eqref{eq:rel_coords_rot} with
feedback~\eqref{eq:CCP} is bounded and the set $\hat W_\delta$ is
globally asymptotically stable. \qed

\section{Simulation Results}\label{sec:simulation}
We consider a group of five robots with the sensor digraph in Figure~\ref{fig:graph}.
The robot masses and inertia matrices are: $m_1=3$ Kg,\ $m_2=3$ Kg,\ $m_3=3.4$ Kg,\ $m_4=3.2$ Kg,\ $m_5=3.2$ Kg and $J_1:= \diag{(0.13,0.13,0.04)}$\,Kg$\cdot$m${}^2$, as in~\cite{Abdessameud}, $J_2=J_1,\ J_3=1.4 J_1,\ J_4=1.2 J_1,\ J_5=1.2 J_1$. 
\begin{table}[t]
\caption{Simulation Initial Conditions}
\centering
 %used for centering table
\begin{tabular}{c c c c}
% centered columns (4 columns)
\hline\hline
Vehicle $i$ & $x_i(0)$ (m) & $v_i(0)$ (m/s) & $R_i(0)$\\% inserts table
\hline
1 & $(0,-10,10)$ & $(0,0,0)$ & side $1$ \\
2 & $(0,10,10)$ & $(0,0,0)$ & side $2$ \\
3 & $(0,0,0)$ & $(0,0,0)$ & down \\
4 & $(-10,0,-10)$ & $(0,0,0)$ & up \\
5 & $(10,0,-10)$ & $(0,0,0)$ & up \\[1ex]
\hline
\end{tabular}
\label{tab:sim}
\end{table}
\begin{table}[t]
\caption{Control Effort}
\centering
 %used for centering table
\begin{tabular}{l c c}
% centered columns (4 columns)
\hline\hline
&Figure~\ref{fig:sim}&Figure~\ref{fig:sim_dist}\\
\hline
$\max_i \sup_t |u_i(t)|$ (N)&20.4&17.21\\
$\max_i \sup_t \|\tau_i(t)\|$ (N$\cdot$m)&15.27&16.47\\
$\max_i \operatorname{rms}(|u_i(t)|)$ (N) &1.72&4.31\\
$\max_i \operatorname{rms}(\|\tau_i(t)\|)$ (N$\cdot$m) &1.43&2.24\\[1ex]
\hline
\end{tabular}
\label{tab:control_effort}
\end{table}
We use the
double-integrator consensus controller of Ren and Atkins~\cite{Ren3}, $\wdes_i(y_i) = \sum_{j=1}^{n}a_{ij}(x_{ij}+\gamma v_{ij})$ where $a_{ij} \geq 0, \ \gamma>0$. It is shown in~\cite{Ren3} that for
sufficiently large $\gamma$ the above controller does indeed achieve
consensus.  We pick $a_{ij}=0.3$ for all $j \in {\cal N}_i$ and
$\gamma=30$. The control gains $k_1$ and $k_2$ in~\eqref{eq:CCP} are
chosen to be $k_1=2$ and $k_2=0.45$. The initial conditions of the robots are shown in
Table~\ref{tab:sim}. The initial attitudes $R_i(0)$ of the robots are: \ up(right), side(ways) $1$, side(ways) $2$ and (upside)down
respectively given by:
\[
\begin{bmatrix}
1 & 0 & 0\\
0 & 1 & 0\\
0 & 0 & 1
\end{bmatrix}, \
\begin{bmatrix}
1 & 0 & 0\\
0 & 0 & -1\\
0 & 1 & 0
\end{bmatrix}, \
\begin{bmatrix}
1 & 0 & 0\\
0 & 0 & 1\\
0 & -1 & 0
\end{bmatrix}, \
\begin{bmatrix}
1 & 0 & 0\\
0 & -1 & 0\\
0 & 0 & -1
\end{bmatrix}.
\]  
Figure~\ref{fig:sim} shows the simulation without the presence of
disturbances while Figure~\ref{fig:sim_dist} shows the simulation when
disturbances are present.  The disturbances are: an additive random
noise with maximum magnitude of $0.25$\,N on the applied force; an additive random
noise with maximum magnitude of $0.25$\,N$\cdot$m on the applied torque; an
additive measurement error for the angular velocity, with maximum
magnitude of $0.25$\,rad/s; an additive random  noise
on the quantity $\wdes_i(y_i^i)$
accounting for errors in measurements of relative displacements and
velocities of the vehicles. The direction of this vector has been
rotated within $0.25$\,rad and the magnitude is scaled between $0.75$ to
$1.25$ times the actual magnitude. The disturbances are updated $10$ times per second. In both cases of Figure~\ref{fig:sim} and Figure~\ref{fig:sim_dist},
the vehicles' positions and velocities converge to a neighborhood of one another. 

In Figure~\ref{fig:sim} the vehicles remain within $0.25$m of one another while in Figure~\ref{fig:sim_dist} the vehicles remain within $1$m of one another at steady state. These neighborhoods can be made even smaller by further increasing the control gains $k_1$ and $k_2$. However, this would result in having higher control inputs. Metrics related to the thrust and torque inputs are presented in Table~\ref{tab:control_effort}. The first two rows show peak control norms and the last two show the root mean square (rms) of the control norms. In these simulations we considered zero gravity, i.e., $g=0$. This was done to improve visibility of the simulation results. In the presence of gravity, the vehicles would still converge to the same neighborhood of one another, however at steady state they would accelerate in the direction of gravity since gravity is not compensated through the control inputs in~\eqref{eq:CCP}. 
\begin{figure}[t]
\begin{center}
\begin{tikzpicture}[->,>=stealth',shorten >=1pt,auto,node distance=1cm,
  thick,main node/.style={circle,draw}]

  \node[main node] (3) {$3$};
  \node[main node] (1) [below of=3] {$1$};
  \node[main node] (2) [above of=3] {$2$};
  \node[main node] (4) [left of=3] {$4$};
  \node[main node] (5) [right of=3] {$5$};
%[every node/.style={font=\sffamily\small}]
    
  \path 
  		(3) edge node[auto] {} (4)
        (3) edge node [left] {} (5)
        (1) edge node [left] {} (3)
        (2) edge node [left] {} (3)
        (4) edge node [left] {} (1)
        (5) edge node [left] {} (2);
\end{tikzpicture}
\end{center}
\caption{Sensor digraph used in the simulation results.}
\label{fig:graph}
\end{figure}
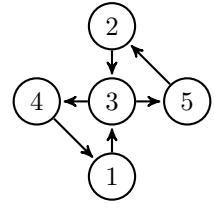
\begin{figure}[t]
\centerline{\includegraphics[width=0.4\textwidth]{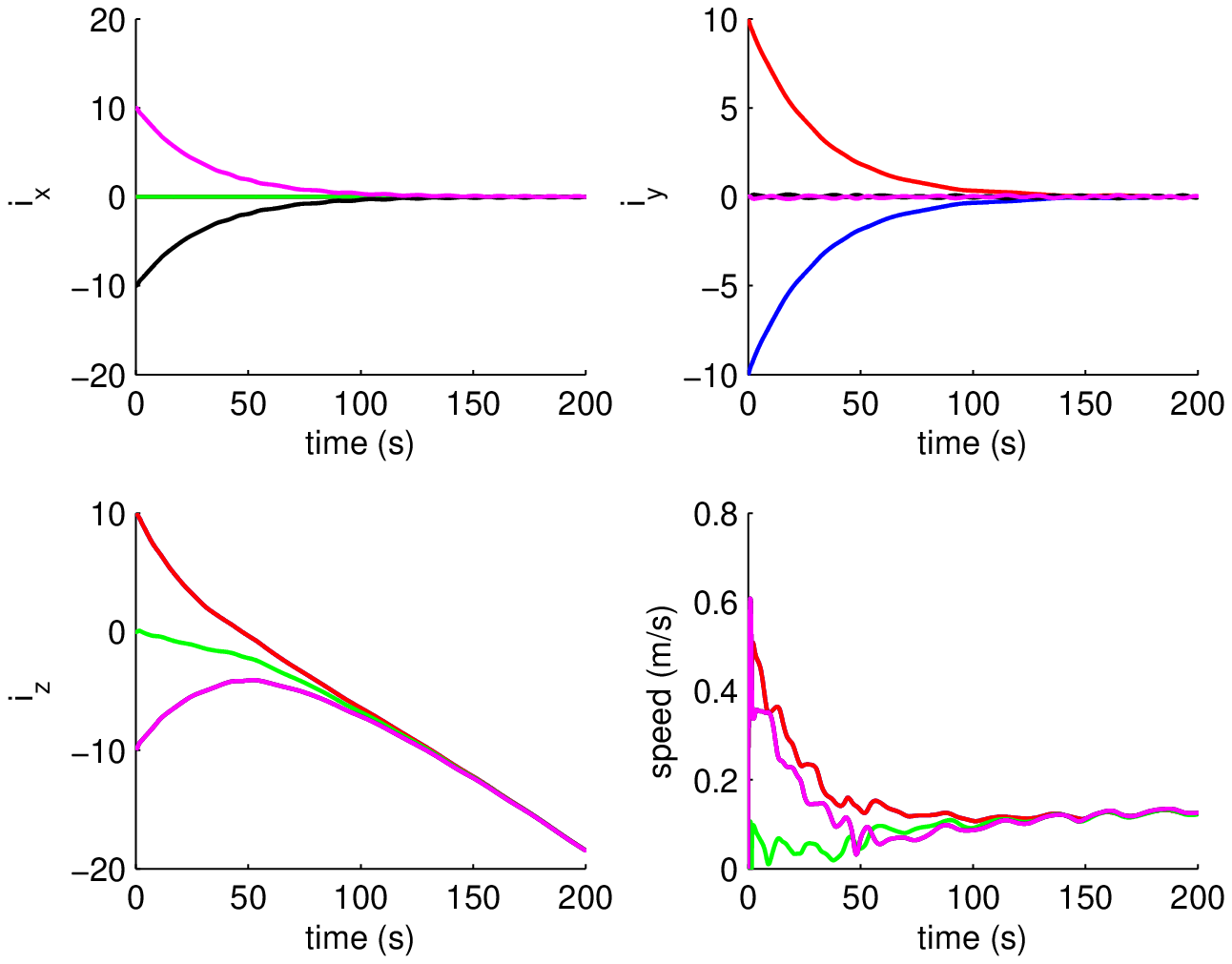}}
\caption{Rendezvous control simulation without the presence of disturbances. At the top-left, top-right and bottom-left: positions of the five robots expressed in the inertial frame ${\cal I}$. At the bottom-right: linear speeds $\|v_{i}\|,\, i=1,\dots,5$.}
\label{fig:sim}
\end{figure}
\begin{figure}[t]
\centerline{\includegraphics[width=0.4\textwidth]{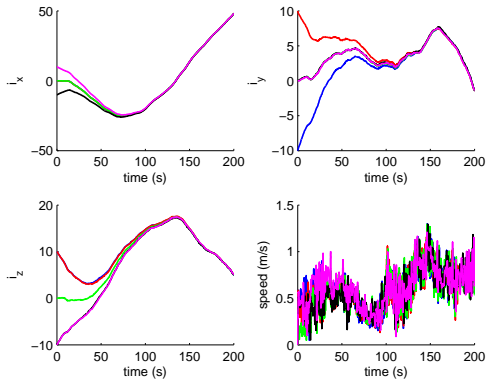}}
\caption{Rendezvous control simulation with the presence of disturbances. At the top-left, top-right and bottom-left: positions of the five robots expressed in the inertial frame ${\cal I}$. At the bottom-right: linear speeds $\|v_{i}\|,\, i=1,\dots,5$.}
\label{fig:sim_dist}
\end{figure}
\section{Conclusions}\label{sec:conclusion}
We have presented the first local and distributed feedback solving the
rendezvous control problem for a class of underactuated robots
modelling vertical take-off and landing (VTOL) vehicles such as
quadrotor helicopters. The main result, Theorem~\ref{thm:2step},
relies on the assumption that the sensor digraph is constant. As we
have discussed in the paper, this assumption is questionable in
practice, 
%but the same is true for the commonly used assumption that
%the sensor digraph is time-varying.  
but a stability analysis in the presence of a state-dependent sensor digraph
is beyond the scope of this paper. We believe that solutions in the literature for consensus of double-integrators with time-dependent sensor digraphs could be extended to rigid bodies using the framework in this paper. However the Lyapunov function used in the analysis would need to be modified extensively. Since this makes the problem even more difficult than it already is, we leave it as a possible future research direction. In this paper we limited ourselves
to the control specification of rendezvous. The proposed control law,
in particular, does not guarantee hovering of the vehicles. While the
robots converge to each other, nothing can be said about the motion of
the ensemble. This cannot be otherwise, for it would be impossible to
solve the rendezvous problem with hovering without additional
sensors. One would need some measurement of the gravity vector, for
example provided by a three-axis accelerometer. The point of view of
these authors is that the proposed solution of the rendezvous problem
will serve as a layer in a hierarchy of higher-level control
specifications such as hovering, formation stabilization, and path
following.

%We have presented a class of feedbacks and a control architecture
%solving the open problem rendezvous control problem for underactuated
%thrust-propelled vehicles using local distributed control. The vehicle class considered has dynamic translational dynamics and kinematic rotational dynamics. The presented control is static and yields global practical stability of the rendezvous configuration. The controller is decomposed into two simpler control stages. In the first stage the vehicles are treated as double-integrator systems . The control revealed several advantages over results in the current literature.   
%We have not explicitly addressed issues of robustness against unmodelled uncertainties and measurement errors but simulation results were presented that include the effect of disturbances on the force and torque inputs. 

\appendix

\subsection{Proof of Lemma~\ref{lem:saturation}}
\noindent Recall the definition of $W(\Chi,R,\omega)$, and assume that $\Chi \neq 0$,
\[
\begin{aligned}
 W =& \alpha \left(\sqrt{V(\Chi)} + \frac 1 2 V(\Chi) \right) +
\sum_{i=1}^n \gdes_i^i(\Chi,R)\cdot e_3 \\
&+ \frac 1 2 (\omega - \OO(\Chi,R)
)\trans \Je (\omega - \OO(\Chi,R)) \\
 \quad =&\sqrt{V(\Chi)} \left( \alpha + \frac{ \sum_{i=1}^n \gdes_i^i(\Chi,R)\cdot
  e_3 }{\sqrt{V(\Chi)}} \right) + \frac \alpha 2 V(\Chi) \\
&  + \frac 1 2 (\omega
- \OO(\Chi,R) )\trans \Je (\omega - \OO(\Chi,R)).
\end{aligned}
\]
Since $\gdes_i^i(\Chi,R)$ is linear with respect to
$\Chi$, we have
\[
\begin{aligned}
W=&\sqrt{V(\Chi)} \left( \alpha + \sum_{i=1}^n
\gdes_i^i\left(\mu(\Chi),R\right)\cdot e_3 \right) + \frac \alpha 2
V(\Chi) \\
&+ \frac 1 2  (\omega - \OO(\Chi,R) )\trans \Je (\omega - \OO(\Chi,R)),
\end{aligned}
\]
where $\mu(\Chi):=\Chi / \sqrt{V(\Chi)}$ is continuous on $\sX \backslash
\{0\}$ and bounded as follows
\[
\| \mu(\Chi) \| = \frac{\| \Chi \|}{\sqrt{V(\Chi)}} = \frac{\|
  \Chi\|}{\sqrt{\Chi\trans P \Chi}} \leq \frac{1}{\sqrt{\lambda_{\min}(P)}}.
\]
Since $\gdes_i^i$ is continuous, $\mu(\Chi)$ is bounded, and $R \in
\sR$, a compact set, it follows that the function $\sum_{i=1}^n \left|
\gdes_i^i\left(\mu(\Chi),R\right)\cdot e_3 \right|$ has a bounded
supremum. Accordingly, let 
\[
\alpha^\star = \sup_{(\Chi,R)\in \sX \backslash \{0\} \times \sR} \sum_{i=1}^n \left|
\gdes_i^i\left(\mu(\Chi),R\right)\cdot e_3 \right|.
\]
For all $\alpha > \alpha^\star$, we have $W(\Chi,R,\omega) \geq \underline{W}(\Chi,R,\omega),$
\[
\begin{aligned}
\underline{W}(\Chi,R,\omega):=& \frac \alpha 2 V(\Chi)\\
&+
\frac 1 2 (\omega - \OO(\Chi,R) )\trans \Je (\omega - \OO(\Chi,R)) \geq 0.
  \end{aligned}
\]
We derived the bound above for $\Chi \neq 0$, but since $\gdes_i^i(0,R)=0$
(by linearity of $\gdes_i^i$ with respect to $\Chi$), the bound also holds
for $\Chi=0$.  The above inequality implies that $W \geq 0$ and
$W^{-1}(0) \subset \underline{W}^{-1}(0)$. But $\underline{W}=0$ if and
only if $V(\Chi)=0$ (i.e., $\Chi=0$) and $\omega = \OO$. Thus $W^{-1}(0)
\subset \Gamma^\star$, proving part (i) of the lemma.

For part (ii), note that for all $c>0$, $W_c \subset \{
\underline{W}\leq c\}$. Since $\underline{W}$ is a positive definite
quadratic form in the variables $(\Chi,\omega - \OO)$, its sublevel
sets are compact in $(\Chi,\omega - \OO)$ coordinates. Thus if
$(\Chi,R,\omega) \in W_c$, $\Chi$ and $\omega - \OO(\Chi,R)$ are
bounded. Since $\OO$ is continuous and $R\in \sR$, a compact set,
$\OO$ is bounded, implying that $\omega $ is also bounded. Therefore
the set $W_c$ is bounded. Continuity of $W$ implies that $W_c$ is compact. This concludes the
proof of part (ii) of the lemma.

For part (iii), let $\varepsilon >0$ be arbitrary. Since
$\underline{W}$ is a positive definite quadratic form in the variables
$(\Chi,\omega - \OO)$, there exists $\delta>0$ such that
$\underline{W}(\Chi,R,\omega) \leq \delta$ implies $\|(\Chi,\omega -
\OO(\Chi,R))\| \leq \varepsilon$.

Furthermore, the inequality $\|(\Chi,\omega -\OO(\Chi,R))\| \leq \varepsilon$ implies that $\|\Chi\| \leq \varepsilon$. Now consider any point $(\Chi,R,\omega)\in \{\underline{W}\leq \delta\}$. We have just seen that this implies that $\|\Chi\| \leq \varepsilon$. It will be shown next that this implies $(\Chi,R,\omega)\in B_\varepsilon(\Gamma^\star)$ and hence $\{\underline{W}\leq \delta\} \subset B_\varepsilon(\Gamma^\star)$. 
\\
Note that $(\Chi,R,\omega)\in \sX \times \sR \times \sO$ lies on the product of metric spaces $\sX$, $\sR$ and $\sO$. Respectively, the metrics are $d_{\sX}$, $d_{\sR}$ and $d_{\sO}$ ($d_{\sX}$ and $d_{\sO}$ are Euclidean metrics). As such, choosing to use the $2$-product metric,
\[
\begin{aligned}
&\|(\Chi,R,\omega)\|_{\Gamma^\star}\\
&=\inf_{(\Chi_0,R_0,\omega_0) \in \Gamma^{\star}}\left(d_{\sX}(\Chi,\Chi_0)^2+d_{\sR}(R,R_0)^2+d_{\sO}(\omega,\omega_0)^2\right)^{\frac{1}{2}}.
\end{aligned}
\]  
Recall that $\Gamma^\star=\{(\Chi,R,\omega) \in \sX \times \sR \times \sO: \Chi=0\}$. As such, the point $(0,R,\omega)$ is contained in the set $\Gamma^{\star}$ and therefore,
\[
\|(\Chi,R,\omega)\|_{\Gamma^\star} \leq \left(d_{\sX}(\Chi,0)^2+d_{\sR}(R,R)^2+d_{\sO}(\omega,\omega)^2\right)^{\frac{1}{2}}
\]
where $d_{\sR}(R,R)$ and $d_{\sO}(\omega,\omega)$ are zero. This yields, $
\|(\Chi,R,\omega)\|_{\Gamma^\star} \leq d_{\sX}(\Chi,0) \leq \|\Chi\| \leq \varepsilon$. This implies that $(\Chi,R,\omega) \in B_\varepsilon(\Gamma^\star)$. Thus, $W_\delta \subset
\{\underline{W}\leq \delta\} \subset B_\varepsilon(\Gamma^\star)$, as
required. This concludes the proof of Lemma~\ref{lem:saturation}.
\qed

\subsection{Proof of Lemma~\ref{lem:Wt}}

We will use a standard result from differential geometry relating the
Lie derivatives of smooth functions along $F$-related vector
fields~\cite[Proposition 8.16]{Lee2}. In our context, recalling that
$\lam = \sqrt V |_{X =\lam \prj}$ and $\hat W = W|_{X = \lam
  \prj}$, the result has the following implication:
\begin{equation}\label{eq:lambda_dot}\dot \lam = \frac{d}{dt} \sqrt{V}\Big|_{\Chi =\lam \theta} \text{ and }
\dot {\hat W}_{\rot} = \frac{d}{dt} W_\rot \Big|_{\Chi = \lam \theta}.
\end{equation}

Rewrite the dynamics of $\Chi$ in~\eqref{eq:rel_coords_trans} as
\[
\begin{aligned}
\dot x_{1j} &= v_{1j} \\
\dot v_{1j} =& \left[ \gdes_j(\Chi) - \gdes_1(\Chi) \right] \\
&+
R_j \left[ (\gdes_j^j(\Chi,R) \cdot e_3) e_3 - \gdes_j^j(\Chi,R) \right]
\\
&+ R_1 \left[ (\gdes_1^1(\Chi,R) \cdot e_3) e_3 -
  \gdes_1^1(\Chi,R) \right].
\end{aligned}
\]
To get the identities above, we added and subtracted
in~\eqref{eq:rel_coords_trans} the ideal force feedbacks
$\wdes_j(y_j)=\gdes_j(\Chi)$ and $\wdes_1(y_1)=\gdes_1(\Chi)$, and we
replaced $u_j$ and $u_1$ in~\eqref{eq:rel_coords_trans} by the
assigned feedbacks in~\eqref{eq:CCP}. Finally, we used the identity
$R_i \gdes_i^i = \gdes_i$.

Taking the time derivative of $\sqrt{V(\Chi)}$ along the above vector
field we get
\[
\begin{aligned}
\frac{d}{dt} \sqrt{V(\Chi)} =& \frac{1}{2 \sqrt{V(\Chi)}} \Bigg[ -\Chi\trans Q \Chi\\
  &+ \sum_{j=2}^n \frac{\partial V}{\partial v_{1j}}  R_j
  \left((\gdes_j^j(\Chi,R) \cdot e_3) e_3 - \gdes_j^j(\Chi,R)\right) \\
&- \sum_{j=2}^n \frac{\partial V}{\partial v_{1j}} R_1
  \left((\gdes_1^1(\Chi,R) \cdot e_3) e_3 \hspace{-.8mm} -\hspace{-.8mm} \gdes_1^1(\Chi,R)\right)
\hspace{-1mm}  \Bigg] \hspace{-.6mm}.
\end{aligned}
\]
The first term in the bracket is the derivative of $V(\Chi)$ along the
nominal vector field~\eqref{eq:consensus_closed_loop}, and $Q=Q\trans$
is a positive definite matrix. Letting $M_2 = \lambda_{\min}(Q) /(2
\lambda_{\max}(P))$ and using the fact that the Euclidean norm is
invariant under rotations, we have
\[
\begin{aligned}
\frac{d}{dt} \sqrt{V(\Chi)} \leq& -M_2 \sqrt{V(\Chi)} + \frac{1}{2
  \sqrt{V(\Chi)}}\cdot \\
&\Bigg[ \sum_{j=2}^n \left\| \frac{\partial V}{\partial
    v_{1j}} \right\| \hspace{-1mm}\left( \| ( \gdes_j^j(\Chi,R) \cdot e_3 ) e_3 -
  \gdes_j^j(\Chi,R) \| \right.\\
  &\left.+ \| ( \gdes_1^1(\Chi,R) \cdot e_3 ) e_3 - \gdes_1^1(\Chi,R)
  \| \right) \Bigg].
\end{aligned}
\]
We claim that $\|( \gdes_i^i(\Chi,R) \cdot e_3 ) e_3 - \gdes_i^i(\Chi,R) \| =
\|\gdes_i^i(\Chi,R) \times e_3 \|$. Indeed, writing $\gdes_i^i = (\gdes_i^i
\cdot e_3) e_3 + \gdes_i^i - (\gdes_i^i \cdot e_3) e_3$, we have $ \gdes_i^i
\times e_3 = ( \gdes_i^i - (\gdes_i^i \cdot e_3) e_3 ) \times e_3$.  Since
the vector $\gdes_i^i - (\gdes_i^i \cdot e_3) e_3$ is perpendicular to
$e_3$, $\|( \gdes_i^i - (\gdes_i^i \cdot e_3) e_3 ) \times e_3 \| =
\|\gdes_i^i - (\gdes_i^i \cdot e_3) e_3 \|$, so that $\|\gdes_i^i \times
e_3\| = \|\gdes_i^i - (\gdes_i^i \cdot e_3) e_3 \|$. This proves the
claim. Using the identity just derived, we get
\[
\begin{aligned}
\frac{d}{dt} \sqrt{V(\Chi)} \leq &-M_2 \sqrt{V(\Chi)}\\
&+ \frac{1}{2
  \sqrt{V(\Chi)}}\Bigg[ \sum_{j=2}^n \left\| \frac{\partial V}{\partial
    v_{1j}} \right\| \left( \| \gdes_j^j(\Chi,R) \times e_3
 \| \right.\\
&\left. + \| \gdes_1^1(\Chi,R) \times e_3 \| \right) \Bigg].
 \end{aligned}
\]

Using~\eqref{eq:lambda_dot}, we get 
\[
\begin{aligned}
\dot \lam \leq -M_2 \lam + \frac{1}{2 \lam} \Bigg[
  \sum_{j=2}^n & \left\| \frac{\partial V}{\partial v_{1j}}(\lam
  \prj) \right\| \left( \| \gdes_j^j(\lam \prj,R) \times e_3 \| \right.\\
  &\left. + \|
  \gdes_1^1(\lam \prj,R) \times e_3 \| \right) \Bigg].
  \end{aligned}
\]
Since the functions $\gdes_i^i$ are linear with respect to their first
argument, and the partial derivatives of the quadratic form $V$ are
linear functions, by the homogeneity  of the norm we have
\[
\begin{aligned}
\dot \lam \leq -M_2 \lam + \frac{\lam}{2} \Bigg[ \sum_{j=2}^n&
  \left\| \frac{\partial V}{\partial v_{1j}}(\prj) \right\| \left( \|
  \gdes_j^j(\prj,R) \times e_3 \| \right.\\
  &\left. + \| \gdes_1^1(\prj,R)
  \times e_3 \| \right) \Bigg].
  \end{aligned}
\]
The functions $\| \partial V/\partial v_{1j} \|$ are continuous. The
variable $\prj$ belongs to $S_1$, a compact set. Therefore $\|
\partial V/\partial v_{1j} \|$ has a maximum,
\[
\begin{aligned}
\dot \lam \leq& -M_2 \lam + \hspace{-1.5mm}\max_{\stackrel{\prj \in S_1
  }{j\in\{2,\ldots,n\}}} \left\| \frac{\partial V}{\partial
  v_{1j}}(\prj) \right\| \frac{\lam}{2} \Bigg[  \sum_{j=2}^n \left(
  \| \gdes_j^j(\prj,R) \times e_3 \| \right) \\
  &+ (n-1) \| \gdes_1^1(\prj,R)
  \times e_3 \| \Bigg] \\
 \leq& - \hspace{-1mm}M_2 \lam +\hspace{-2mm} \max_{\stackrel{\prj \in S_1
  }{j\in\{2,\ldots,n\}}} \hspace{-1mm}\left\| \frac{\partial V}{\partial
  v_{1j}}(\prj) \right\| \frac{\lam}{2} (n-1) \hspace{-1mm}\sum_{j=1}^n \|
\gdes_j^j(\prj,R) \hspace{-1mm}\times\hspace{-1mm} e_3 \|.
\end{aligned}
\]
Letting $M_1:=\max_{\stackrel{\prj \in S_1 }{j\in\{2,\ldots,n\}}} \left\|
\partial V / \partial v_{1j} \right\| (n-1) /2,$ we get the first inequality in~\eqref{eqn:inequalities}.

%%%%%%%%%%%%%%%%%%%%%
%%%%%%%%%%%%%%%%%%%%%
%%%%%%%%%%%%%%%%%%%%%

We now turn to the second inequality
in~\eqref{eqn:inequalities}. Recall the definition of $W_\rot$,
\[
\begin{aligned}
W_\rot(\Chi,R,\omega) =& \sum_{i=1}^n \gdes_i^i(\Chi,R) \cdot e_3 \\
&+ \frac 1 2
(\omega - \OO(\Chi,R) ) \trans \Je (\omega - \OO(\Chi,R)).
\end{aligned}
\]
The time derivative of $W_\rot$  along the vector field
in~\eqref{eq:rel_coords_trans}-\eqref{eq:rel_coords_rot} is
\[
\begin{aligned}
\dot W_\rot = \sum_{i=1}^n \Bigg[ &\left( \frac{d}{dt} \gdes_i^i \right)
  \cdot e_3 + (\omega_i^i - \OO_i^i(\Chi,R))\cdot \\
  &  \left(\tau_i -
  \omega_i^i \times J_i \omega_i^i - J_i \left( \frac{d}{dt} \OO_i^i
  \right) \right) \Bigg].
  \end{aligned}
\]
To express $(d / dt) \gdes_i^i$, recall that $\gdes_i^i(\Chi,R) =
R_i^{-1} \wdes_i(h_i(\Chi))$. Then,
\[
\frac d {dt} \gdes_i^i = \left( \frac d {dt} R_i^{-1} \right)
\wdes_i(h_i(\Chi)) + R_i^{-1} \frac d {dt} \left( \wdes_i(h_i(\Chi))
\right).
\]
The function $\wdes_i(h_i(\Chi))$ is linear. Its derivative along the
vector field~\eqref{eq:rel_coords_trans}-\eqref{eq:rel_coords_rot}
with feedback~\eqref{eq:CCP} is a function of $(\Chi,R)$ which is linear
with respect to $\Chi$ because $u_i = - \gdes_i^i(\Chi,R) \cdot e_3$ is
such. We will denote it $\hdes_i(\Chi,R)$, $\hdes_i(\Chi,R) :=(d /dt)
\wdes_i(h_i(\Chi))$. Consistently with our notational convention in
Table~\ref{table:abs_rel}, we will let $\hdes_i^i(\Chi,R):= R_i^{-1}
\hdes_i(\Chi,R)$. The function $\hdes_i^i(\Chi,R)$ is linear with respect to
$\Chi$. Returning to the derivative of $\gdes_i^i$, we have
\[
\begin{aligned}
\frac d {dt} \gdes_i^i &= - (\omega_i^i)^\times R_i^{-1} \wdes_i(h_i(\Chi))
+ R_i^{-1} \hdes_i(\Chi,R) \\
&= - \omega_i^i \times \gdes_i^i(\Chi,R) + \hdes_i^i(\Chi,R).
\end{aligned}
\]
Similarly, since $\OO_i^i(\Chi,R)=k_1 (\gdes_i^i(\Chi,R) \times e_3)$, we have $\frac{d}{dt} \OO_i^i = k_1 \left(- \omega_i^i \times \gdes_i^i +
\hdes_i^i \right) \times e_3$. Substituting the above identities in the expression for $\dot W_\rot$
and since $\tau_i = \omega_i^i \times J_i \omega_i^i - k_1 J_i
((\omega_i^i \times \gdes_i^i) \times e_3) - k_1^2 k_2 (\omega_i^i -
\OO_i^i)$, we get
\[
\begin{aligned}
\dot W_\rot =& \sum_{i=1}^n \left[ -(\omega_i^i \times \gdes_i^i)\cdot
  e_3 + \hdes_i^i \cdot e_3 \right.\\
   &\left. -k_1 (\omega_i^i - \OO_i^i) \cdot  J_i (
  \hdes_i^i \times e_3) - k_1^2 k_2 \| \omega_i^i - \OO_i^i\|^2
  \right].
\end{aligned}
\]
Using the property of the triple product that $(\omega_i^i \times
\gdes_i^i) \cdot e_3 = (\gdes_i^i \times e_3) \cdot \omega_i^i$, we obtain
\[
\begin{aligned}
\dot W_\rot = &\sum_{i=1}^n \left[ - (\gdes_i^i \times e_3) \cdot
  \omega_i^i+ \hdes_i^i \cdot e_3 \right.\\ 
    &\left. -k_1 (\omega_i^i - \OO_i^i) \cdot J_i
  ( \hdes_i^i \times e_3)   
 - k_1^2 k_2 \| \omega_i^i - \OO_i^i\|^2
  \right].
\end{aligned}
\]
Adding and subtracting the term $(\gdes_i^i \times e_3) \cdot \OO_i^i$
and collecting the term $\omega_i^i - \OO_i^i$, we have
\[
\begin{aligned}
\dot W_\rot =& \sum_{i=1}^n \big[ - (\gdes_i^i \times e_3) \cdot \OO_i^i
  + \hdes_i^i \cdot e_3 - ( (\gdes_i^i \times e_3) \\
&  + k_1 J_i ( \hdes_i^i
  \times e_3) ) \cdot (\omega_i^i - \OO_i^i) 
 - k_1^2 k_2 \| \omega_i^i - \OO_i^i\|^2 \big].
\end{aligned}
\]
Substituting in the first term inside the bracket $\OO_i^i=-k_1
(\gdes_i^i \times e_3)$, taking norms, and using the fact that $k_1 \geq
1$, we arrive at the inequality
\[
\begin{aligned}
\dot W_\rot  \leq &\sum_{i=1}^n \big[ - k_1 \| \gdes_i^i \times e_3\|^2
  + \| \hdes_i^i \cdot e_3\| \\
  &+ k_1 \kdes_i(\Chi,R) \| \omega_i^i -
  \OO_i^i\| - k_1^2 k_2 \| \omega_i^i - \OO_i^i\|^2 \big],
  \end{aligned}
\]
where $\kdes_i(\Chi,R):=\|\gdes_i^i(\Chi,R) \times e_3\| + \|J_i (
\hdes_i^i(\Chi,R) \times e_3)\|$. Note that $\kdes_i(\Chi,R)$ is
homogeneous with respect to $\Chi$ because $\gdes_i^i$ and $\hdes_i$
are linear with respect to $\Chi$ and the norm is a homogeneous
function.

Splitting the term $- k_1^2 k_2 \| \omega_i^i - \OO_i^i\|^2$ into two
parts and noticing that the function $k_1 \kdes_i(\Chi,R) \| \omega_i^i -
\OO_i^i\| - (k_1^2 k_2/2) \| \omega_i^i - \OO_i^i\|^2$ is quadratic in
the variable $\| \omega_i^i - \OO_i^i\|$ with maximum $\kdes_i^2(\Chi,R) /
(2 k_2)$, we get
\[
\begin{aligned}
\dot W_\rot  \leq &\sum_{i=1}^n \bigg[ - k_1 \| \gdes_i^i \times e_3\|^2
  + \| \hdes_i^i \cdot e_3\| - \frac{k_1^2 k_2}{2} \| \omega_i^i -
  \OO_i^i\|^2 \\
  &+ \frac{ \kdes_i^2(\Chi,R) }{2 k_2} \bigg].
  \end{aligned}
\]
Now using~\eqref{eq:lambda_dot} we get
\[
\begin{aligned}
\dot {\hat W}_\rot \leq & \sum_{i=1}^n \bigg[ - k_1 \| \gdes_i^i(\lam
  \prj,R) \times e_3\|^2 + \| \hdes_i^i(\lam \prj,R) \cdot e_3\|\\
  & -
  \frac{k_1^2 k_2}{2} \| \omega_i^i - \OO_i^i\|^2 + \frac{
    \kdes_i^2(\lam \prj,R) }{2 k_2} \bigg].
      \end{aligned}
\]
Using the homogeneity with respect to $\Chi$ of $\|\gdes_i^i \times e_3\|$,
$\| \hdes_i^i \cdot e_3\|$, and $\kdes_i$, we get
\[
\begin{aligned}
\dot {\hat W}_\rot \leq & \sum_{i=1}^n \bigg[ - k_1 \lam^2 \|
  \gdes_i^i(\prj,R) \times e_3\|^2 + \lam \| \hdes_i^i(\prj,R) \cdot
  e_3\| \\
  &- \frac{k_1^2 k_2}{2} \| \omega_i^i - \OO_i^i\|^2 + \lam^2
  \frac{ \kdes_i^2(\prj,R) }{2 k_2} \bigg].
    \end{aligned}
\]
Since $\| \hdes_i^i(\prj,R) \cdot e_3\|$ and $\kdes_i^2(\prj,R)$ are
continuous functions over the compact set $S_1 \times \sR$, they each
have a maximum. Letting $
M_3 = n \cdot \max_{\stackrel{(\theta,R)\in S_1 \times
    \sR}{i\in\{1,\ldots,n\}}}\hspace{-1mm} \left( \| \hdes_i^i(\prj,R)
\cdot e_3\| \right)$,
$M_4 =\max_{\stackrel{(\theta,R)\in S_1 \times
    \sR}{i\in\{1,\ldots,n\}}} \hspace{-1mm}\left(\hspace{-0.6mm}
\frac{\kdes_i^2(\prj,R)}{2} \hspace{-0.8mm}\right),
$ we conclude that
\[
\begin{aligned}
\dot {\hat W}_\rot \leq &\lam^2 \sum_{i=1}^n \left[ -k_1 \|
  \gdes_i^i(\prj,R) \times e_3\|^2 + \frac{M_4}{k_2} \right] + \lam
M_3 \\
&-  \frac{k_1^2 k_2}{2} \sum_{i=1}^n\| \omega_i^i - \OO_i^i\|^2,
  \end{aligned}
\]
as required. This concludes the proof of Lemma~\ref{lem:Wt}.
\qed
\bibliographystyle{IEEEtran} \bibliography{sub}
\begin{IEEEbiography}[{\includegraphics[width=1.1in,height=1.25in,clip,keepaspectratio]{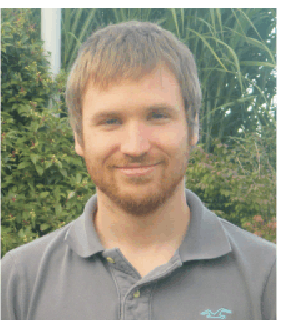}}]{Ashton Roza} received the Bachelor of Applied Science (BASc) degree in electrical engineering in 2010 from the University of Waterloo, Canada and the Master of Applied Science (MASc) from the University of Toronto, Canada, in 2012. He is currently a PhD student in the Edward S. Rogers Sr. Department of Electrical and Computer Engineering, University of Toronto, Canada. His research focuses on the application of nonlinear control and geometric methods to coordination in networks of autonomous robots.
\end{IEEEbiography}
\begin{IEEEbiography}[{\includegraphics[width=1.1in,height=1.25in,clip,keepaspectratio]{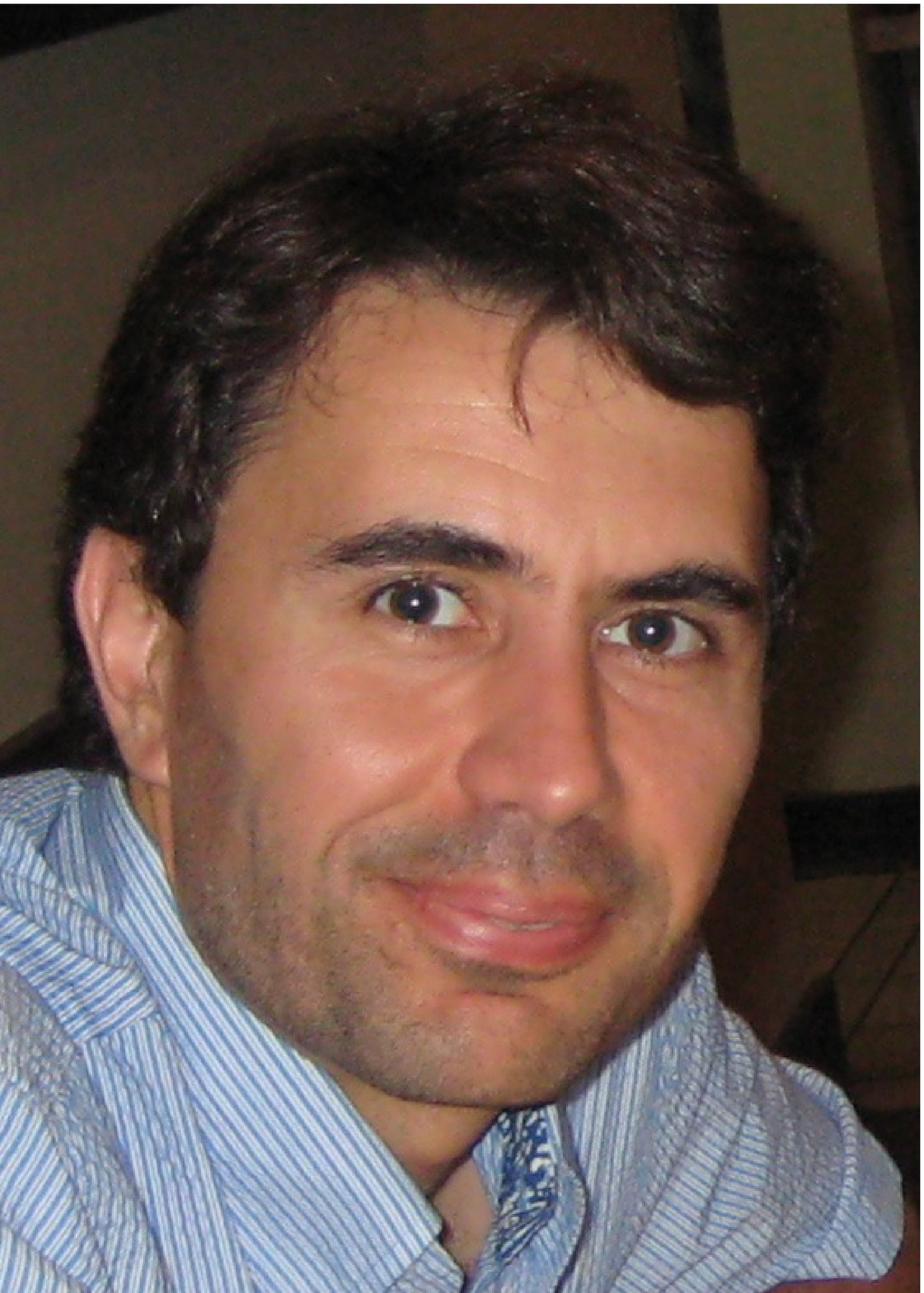}}]{Manfredi Maggiore} was born in Genoa, Italy. He received the "Laurea"
degree in Electronic Engineering in 1996 from the University of Genoa
and the PhD degree in Electrical Engineering from the Ohio State
University, USA, in 2000. Since 2000 he has been with the Edward
S. Rogers Sr. Department of Electrical and Computer Engineering,
University of Toronto, Canada, where he is currently Professor. He has
been a visiting Professor at the University of Bologna (2007-2008),
and the Laboratoire des Signaux et Syst\`emes, Ecole CentraleSup\'elec
(2015-2016). His research focuses on mathematical nonlinear control,
and relies on methods from dynamical systems theory and differential
geometry.
\end{IEEEbiography}
\begin{IEEEbiography}[{\includegraphics[width=1.1in,height=1.25in,clip,keepaspectratio]{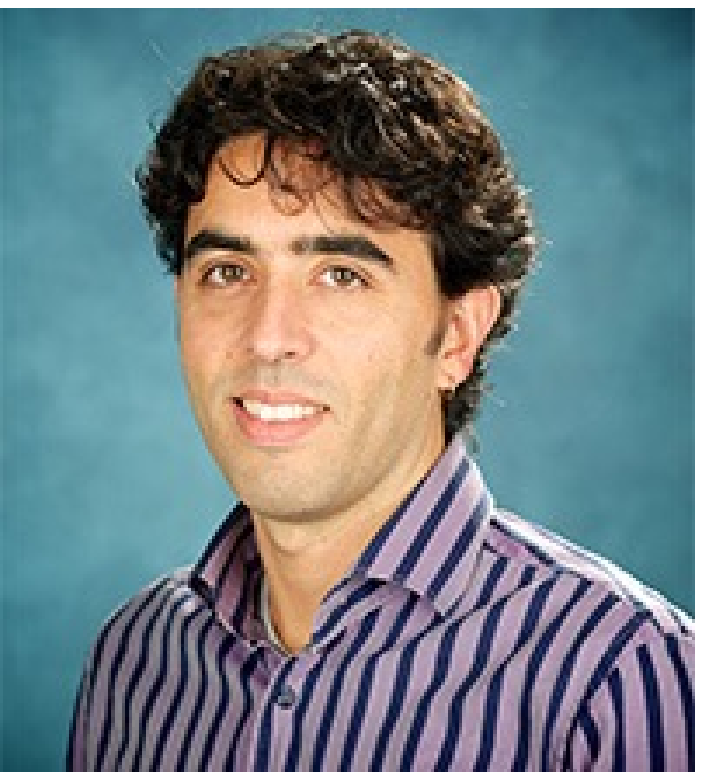}}]{Luca Scardovi}(M '06)
received the Laurea degree and Ph.D. degree in Electronic and Computer Engineering from the University of Genoa, Italy, in 2001 and 2005 respectively. In 2005 he was an Adjunct Professor at the University of Salento, Lecce, Italy. He held research associate positions at the Department of Electrical Engineering and Computer Science at the University of Li\`ege, Belgium (2005-2007) and at the Department of Mechanical and Aerospace Engineering at Princeton University (2007-2009). From 2009 to 2011 he was an Assistant Professor at the Technische Universit\"at M\"unchen (TUM), Munich, Germany. He is currently an Assistant Professor in the Edward S. Rogers Sr. Department of Electrical and Computer Engineering at the University of Toronto. His research interests focus on dynamical systems with special emphasis in the analysis and control of emergent dynamics in networked dynamical systems. 
He is an associate editor for the  {\em IEEE Control Systems Society Conference Editorial Board} and for the journal {\em Systems \& Control Letters}.
\end{IEEEbiography}
\end{document}